\documentclass{gtart}

\def\ifplaintex{\expandafter\ifx\csname documentclass\endcsname\relax}

\expandafter\ifx\csname epsfbox\endcsname\relax\input epsf\fi

\ifplaintex 
\else
\fi

\def\gt{{\mathsurround=0pt\it $\cal G\mskip-2mu$eometry \&\ 
$\cal T\!\!$opology}}        

\def\gtp{{\mathsurround=0pt\it $\cal G\mskip-2mu$eometry \&\ 
$\cal T\!\!$opology $\cal P\!$ublications}}  

\def\lognumber#1{\def\thelognumber{#1}}
\def\volumenumber#1{\def\thevolumenumber{#1}}
\def\papernumber#1{\def\thepapernumber{#1}}
\def\volumeyear#1{\def\thevolumeyear{#1}}

\def\pagenumbers#1#2{\def\startpage{#1}\def\finishpage{#2}}
\def\published#1{\def\publishdate{#1}}
\def\proposed#1{\def\theproposer{#1}}
\def\seconded#1{\def\theseconders{#1}}
\def\received#1{\def\receiveddate{#1}}
\def\revised#1{\def\reviseddate{#1}}
\def\accepted#1{\def\accepteddate{#1}}
\def\asciititle#1{\def\theasciititle{#1}}

\def\asciiaddress#1{\def\theasciiaddress{#1}}
\def\asciiemail#1{\def\theasciiemail{#1}}
\long\def\asciiabstract#1{\long\def\theasciiabstract{#1}}
\def\asciikeywords#1{\def\theasciikeywords{#1}}

\let\\\par\let\thevolumenumber\relax\let\thepapernumber\relax
\let\thevolumeyear\relax\let\thesamplenumber\relax\let\startpage\relax
\let\finishpage\relax\let\publishdate\relax\let\receiveddate\relax
\let\reviseddate\relax\let\accepteddate\relax\let\theasciititle\relax
\let\theasciiauthors\relax\let\theasciiaddress\relax
\let\theasciiabstract\relax\let\theasciikeywords\relax
\let\theasciiemail\relax\let\theshortauthors\relax\let\theshorttitle\relax

\long\def\maketitlep{   

\count0=\startpage

\gt\hfill      
\hbox to 77pt{\vbox to 0pt{\vglue -15pt\epsfbox{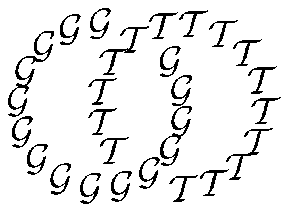}\vss}\hss}
\break
{\small\ifx\thesamplenumber\relax 
Volume \else Sample
\fi\thevolumenumber\ (\thevolumeyear)
\startpage--\finishpage\nl
Published: \publishdate}
\vglue 0.5truein plus 0.4fil minus 0.1truein

{\parskip=0pt\leftskip 0pt plus 1fil\def\\{\par\smallskip}{\ifplaintex\large
\else\Large\fi\bf\thetitle}\par\medskip}   

\vglue 0pt plus 0.1fil 

{\parskip=0pt\leftskip 0pt plus 1fil\def\\{\par}{\sc\theauthors}
\par\medskip}

\vglue 0pt plus 0.1fil 

{\small\parskip=0pt\let\newline\\
{\leftskip 0pt plus 1fil\def\\{\par}{\sl\theaddress}\par}
\expandafter\ifx\theemail\relax    
\relax\else\vglue 5pt plus 0.02fil minus 2pt\def\\{\stdspace{\rm 
and}\stdspace} 
\cl{Email:\stdspace\tt\theemail}\fi
\ifx\theurl\relax                  
\relax\else\vglue 5pt plus 0.02fil minus 2pt\def\\{\stdspace{\rm 
and}\stdspace}
\cl{URL:\stdspace\tt\theurl}\fi\par}

\vglue 7pt plus 0.3fil minus 3pt

{\bf Abstract}
\vglue 5pt plus 0.1fil minus 2pt

\theabstract

\vglue 7pt plus 0.3fil minus 3pt

{\bf AMS Classification numbers}\quad Primary:\quad \theprimaryclass

Secondary:\quad \thesecondaryclass

\vglue 5pt plus 0.3fil minus 2pt

{\bf Keywords:}\quad \thekeywords

\vglue 10pt plus 0.5fil minus 5pt

{\small  Proposed: \theproposer\hfill Received: \receiveddate\nl
Seconded: \theseconders\hfill 
\ifx\reviseddate\relax                         
Accepted: \accepteddate                        
\else
Revised: \reviseddate                          
\fi}
\eject
}       

\let\maketitlepage\maketitlep
\let\maketitle\maketitlepage

\font\phead=cmsl9 scaled 950
\font=cmsl9 scaled 1050
\font\pnum=cmbx10 scaled 913
\font\lnum=cmbx10 
\font\pfoot=cmsl9 scaled 950
\font=cmsl9 scaled 1050
\ifplaintex
\headline{\vbox to 0pt{\vskip -4.5mm\line{\small\phead\ifnum
\count0=\startpage ISSN 1364-0380 (on line)
1465-3060 (printed) \hfill {\pnum\folio}\else\ifodd\count0\def\\{ }%
\ifx\theshorttitle\relax\thetitle\else\theshorttitle\fi\hfill{\pnum\folio}
\else\def\\{ and }{\pnum\folio}\hfill\ifx\theshortauthors\relax\theauthors
\else\theshortauthors\fi\fi\fi}\vss}}
\footline{\vbox to 0pt{\vglue 0mm\line{\small\pfoot\ifnum\count0=\startpage
\copyright\ \gtp\hfill\else
\gt, Volume \thevolumenumber\ (\thevolumeyear)\hfill\fi}\vss
}}
\else
\makeatletter
\def\@oddhead{{\small\lhead\ifnum\count0=\startpage ISSN 1364-0380 (on line)
1465-3060 (printed) \hfill {\lnum\number\count0}\else\ifodd\count0
\def\\{ }\ifx\theshorttitle\relax \thetitle \else\theshorttitle\fi\hfill
{\lnum\number\count0}\else\def\\{ and }{\lnum\number\count0}
\hfill\ifx\theshortauthors\relax 
\theauthors\else\theshortauthors\fi\fi\fi}}\def\@evenhead{@oddhead}
\def\@oddfoot{\small\lfoot\ifnum\count0=\startpage\copyright\ \gtp\hfill\else
\gt, Volume \thevolumenumber\ (\thevolumeyear)\hfill\fi}
\def\@evenfoot{@oddfoot}
\makeatother
\fi

\newwrite\gtoutfile
\long\gdef\makeheadfile{  
{\def\\{, }\def\s{ }
\immediate\openout\gtoutfile head.xxx
\immediate\write\gtoutfile{To: math@arxiv.org}
\immediate\write\gtoutfile{Subject: put OR rep NNNNN:pppp}
\immediate\write\gtoutfile{--text follows this line--}
\immediate\write\gtoutfile{Proxy-for: \ifx\theasciiauthors\relax
\theauthors\else\theasciiauthors\fi\s<\ifx\theasciiemail\relax\theemail\else\theasciiemail\fi>}
\immediate\write\gtoutfile{\noexpand\\}
\immediate\write\gtoutfile{Authors: \ifx\theasciiauthors\relax
\theauthors\else\theasciiauthors\fi}
{\def\\{ }\immediate\write\gtoutfile{Title: \ifx\theasciititle\relax
\thetitle\else\theasciititle\fi}}
\immediate\write\gtoutfile{Subj-class: GT or GR or SG or ...}
\immediate\write\gtoutfile{MSC-class: \theprimaryclass\ifx\thesecondaryclass\relax\else, \thesecondaryclass\fi}
\immediate\write\gtoutfile{Journal-ref: Geom. Topol. \thevolumenumber\s
(\thevolumeyear) \startpage-\finishpage}
\immediate\write\gtoutfile{Comments: Published in Geometry and Topology at}
\immediate\write\gtoutfile{    http://www.maths.warwick.ac.uk/gt/GTVol\thevolumenumber/paper\thepapernumber.abs.html}
\immediate\write\gtoutfile{\noexpand\\}
\immediate\write\gtoutfile{}
\ifx\theasciiabstract\relax
\immediate\write\gtoutfile{\theabstract}\else
\immediate\write\gtoutfile{\theasciiabstract}\fi
\immediate\write\gtoutfile{}
\immediate\write\gtoutfile{\noexpand\\}
\immediate\write\gtoutfile{}
\immediate\closeout\gtoutfile}}  

\def\maketitlepage{\maketitlep\makeheadfile}
\let\maketitle\maketitlepage

\lognumber{121}
\volumenumber{5}\papernumber{18}\volumeyear{2001}
\pagenumbers{551}{578}
\received{31 May 2000}
\revised{2 May 2001}
\accepted{9 June 01}
\proposed{Robion Kirby}
\seconded{Joan Birman, Cameron Gordon}
\published{17 June 2001}
\usepackage{pb-diagram,lamsarrow,pb-lams,amsmath,amssymb,graphicx}

\theoremstyle{plain}

\newtheorem{theorem}{Theorem}
\newtheorem{proposition}{Proposition}[section]
\newtheorem{lemma}[proposition]{Lemma}

\theoremstyle{definition}
\newtheorem{definition}[proposition]{Definition}

\theoremstyle{remark}
\newtheorem{example}[proposition]{Example}

\newtheorem{remark}[proposition]{Remark}

\newcommand{\psdraw}[2]
{\centerline{\includegraphics[width=#2]{#1.eps}}}
\newcommand{\ppsdraw}[2]
{\includegraphics[width=#2]{#1.eps}}

\newcommand{\noi}{\noindent}
\newcommand{\out}{\operatorname{Out}}

\def\lbl#1{\label{#1}}

\def\BZ{\mathbb Z}
\def\BQ{\mathbb Q}

\def\BC{\mathbb C}

\def\A{\mathcal A}

\def\K{\mathcal K}
\def\KY{\mathcal K^{\mathsf{Y}}}

\def\G{\mathcal G}

\def\F{\mathcal F}

\def\l{\lambda}

\def\ga{\gamma}

\def\ihs{integral homology 3--sphere}

\def\fti{finite type
invariant}

\def\Fas#1{\mathcal F_{#1}(N)}

\def\FY#1{\mathcal F^{\mathsf{Y}}_{#1}(N)}

\def\Gas#1{\mathcal G_{#1}(N)}

\def\GY#1{\mathcal G^{\mathsf{Y}}_{#1}(N)}

\def\Las{\mathcal L^{\mathrm{as}}}
\def\Lns{\mathcal L^{\mathrm{ns}}}

\def\s{$ \spadesuit $}

\def\edge{\mathrm{Edge}}
\def\Zpi{\BZ\pi}
\def\simk{\sim_\kappa}
\def\simkp{\sim_{\kappa'}}
\def\Deg{\mathrm{Deg}}

\def\i{^{-1}}
\def\Z{\mathbb Z}

\def\a{\alpha}

\def\bd{\partial}
\def\g{\gamma}

\def\e{\epsilon}

\def\b{\beta}

\def\s{\sigma}

\def\sub{\subset}
\def\lk{{\text{lk}}}

\def\H{\mathcal H}

\renewcommand{\ker}{\operatorname{Ker}}

\newcommand{\con}{\equiv}

\def\sminus{\smallsetminus}
\def\ti{\widetilde}

\def\AS{\mathrm{AS}}
\def\IHX{\mathrm{IHX}}

\def\ygraph{$\mathrm{Y}$--graph}
\def\ylink{$\mathrm{Y}$--link}
\def\clover{clover}

\def\Hom{\mathrm{Hom}}
\newcommand{\oZ}{\otimes\mathbb Z[1/2]}

\renewcommand{\ker}{\operatorname{Ker}}

\def\a{\alpha}
\def\b{\beta}

\def\e{\epsilon}

\def\g{\gamma}
\def\l{\lambda}

\def\Z{\mathbb Z}

\def\H{\mathcal H}

\def\F{\mathcal F}
\def\K{\mathcal K}
\def\A{\mathcal A}

\def\L{{\mathcal L}}

\def\G{\mathcal G}

\def\bd{\partial}

\def\i{^{-1}}

\def\iso{\cong}

\def\con{\equiv}

\def\wt{\widetilde}

\begin{document}

\title
{Homology surgery and invariants of 3--manifolds}
\asciititle
{Homology surgery and invariants of 3-manifolds}

\author{Stavros Garoufalidis\\Jerome Levine}
\address{School of Mathematics\\Georgia Institute of
Technology\\Atlanta, 
GA 30332-0160, USA}
\email{stavros@math.gatech.edu}
\secondaddress{Department of Mathematics\\Brandeis
University\\Waltham, MA 02254-9110, USA}
\secondemail{levine@brandeis.edu}

\asciiaddress{School of Mathematics\\Georgia Institute of
Technology\\Atlanta, 
GA 30332-0160, USA\\Department of Mathematics\\Brandeis
University\\Waltham, MA 02254-9110, USA}
\asciiemail{stavros@math.gatech.edu, levine@brandeis.edu}

\begin{abstract} 
We introduce a homology surgery problem
in dimension 3 which has the property that the vanishing of its
algebraic obstruction leads to a canonical class of $\pi$--algebraically-split
links in 3--manifolds with fundamental group $\pi$.
Using this class of links, we define a theory of \fti s of
3--manifolds in such a way that invariants of degree $0$ are precisely
those of conventional algebraic topology and surgery theory. When \fti
s are reformulated in terms of {\em clovers}, we deduce upper bounds
for the number of invariants in terms of $\pi$--decorated trivalent
graphs. We also consider an associated notion of surgery equivalence of
$\pi$--algebraically split links and prove a classification theorem using a 
generalization of Milnor's $\bar\mu$--invariants to this class of links.
\end{abstract}

\asciiabstract{We introduce a homology surgery problem in dimension 3
which has the property that the vanishing of its algebraic obstruction
leads to a canonical class of \pi-algebraically-split links in
3-manifolds with fundamental group \pi.  Using this class of links, we
define a theory of finite type invariants of 3-manifolds in such a way
that invariants of degree 0 are precisely those of conventional
algebraic topology and surgery theory. When finite type invariants are
reformulated in terms of clovers, we deduce upper bounds for the
number of invariants in terms of \pi-decorated trivalent graphs. We
also consider an associated notion of surgery equivalence of
\pi_algebraically split links and prove a classification theorem using
a generalization of Milnor's \mu-invariants to this class of links.}

\asciikeywords{Homology surgery, finte type invariants, 3-manifolds,
clovers}

\maketitlepage

\section{Introduction}
\lbl{sec.intro}

In this paper we take a new approach to the role of finite-type invariants in
$3$--manifold topology. Our approach is to subdivide the collection of
$3$--manifold invariants into three increasingly delicate classes.
\begin{enumerate}
\item The invariants of {\em classical} algebraic topology, which we take to
mean invariants of {\em homology type} in the strongest sense.
\item Surgery-theoretic invariants.
\item Invariants which we consider to be of {\em finite-type}.
\end{enumerate}
In this spirit we start with a base manifold $N$ and then consider all manifolds
which are homology equivalent to it, over $\pi_1 (N)$. This puts us in a
(homology) surgery-theoretic framework with a resulting Witt-type invariant. The
vanishing of this invariant then restricts us to a class of manifolds which can
be constructed from $N$ by surgery on framed links in $N$ which are
algebraically split in a suitable sense. We can, within this class, define a
notion of finite-type invariant analogous to earlier
notions which were considered, most effectively, for the class of homology
spheres
(corresponding to $N=S^3$).

A natural question to ask is whether this class of manifolds, and the notion of
finite type, depends only on the homology equivalence class of $N$. We make
some progress toward an affirmative answer.

We then give a reformulation of this finite-type theory in terms of what has
been recently called  {\em \ygraph s} (see \cite{Gu}) or
{\em claspers} (see \cite{H}) or {\em \clover s} (see \cite{GGP}). 
This will enable us to deduce upper bounds
for the number of invariants of a given degree from the number of
$\pi$--decorated trivalent
graphs of that degree.

Finally we consider a notion of surgery equivalence for {\em algebraically
split} links in a general $3$--manifold, which is closedly related to our
finite-type theory, generalizing the relation between classical surgery
equivalence in $S^3$ and finite-type theory in homology spheres, as explained 
in \cite{GL}. We then define  $\Z [\pi\times\pi ]$--valued triple Milnor 
invariants and show that they classify surgery equivalence, generalizing 
\cite{Le}. We also give a direct proof that concordant links are surgery
equivalent.

\section{A surgery problem}
\lbl{sec.surgery}

Throughout this paper, all manifolds will be smooth and oriented and all 
maps will be orientation preserving. 
Let $N$ be a closed $3$--manifold. Consider {\em degree} $1$ maps 
$f\co M\to N$, where $M$ is also a closed oriented $3$--manifold. Then the 
induced $f_{\ast}\co \pi_1 (M)\to\pi_1(N )=\pi$ is onto and we can
consider the 
induced homomorphism $f_{\ast}\co H_{\ast}(\ti M)\to H_{\ast}(\ti N )$, 
where $\ti N ,\ti M$ indicates the $\pi$--coverings. We will say that $f$
is 
a $\Zpi$--{\em homology equivalence} if $f_{\ast}$ (on $H_* (\ti M)$) is an
isomorphism and $f$ is degree $1$.
Since $H_1 (\ti N)=0$ and $\pi_1 (\ti M)\iso\ker f_*$ (on $\pi_1 (M)$), it
follows
by Poincar\'e duality (see \cite[Lemma 2.2]{Wa}) that this is equivalent to the
condition that $\ker f_{\ast}$
is a 
{\em perfect} subgroup. Given another $\Zpi$--homology equivalence 
$f'\co M'\to N$, we say they are  {\em diffeomorphically equivalent} iff
there 
exists a diffeomorphism $g\co M\to M'$ such that $f$ is homotopic to
$f'\circ g$.
Let $\H(N)$ denote the {\em structure set} of {\em diffeomorphism 
equivalence classes} and let $\H_0(N)$ (resp. $\H_0^s (N)$) denote the 
set of (simple) $\Zpi$--homology bordism classes of (simple) $\Zpi$--homology 
equivalences $f\co M\to N$.

Our goal in this section is to define a {\em (homology)  surgery obstruction
map} $\Phi$ and its relatives $\Phi_0$ and $\Phi_0^s$ which fit in the
following commutative diagram:

\begin{theorem}
\lbl{thm.cd}
$$
\begin{diagram}
\node{\H_0^s (N)}\arrow{e}\arrow{s,l}{\Phi_0^s}\node{\H_0 (N
)}\arrow{s,l,A}{\Phi_0}\node{\H (N)}\arrow{w,A}\arrow{s,l,A}{\Phi}\\
\node{\wt W_s (\pi )}\arrow{e}\node{\wt W(\pi )}\node{B (\pi
)}\arrow{w,A}
\end{diagram}
$$
\end{theorem}

\subsection{Algebraic preliminaries}
\lbl{sub.alg}

In this section we define the semigroups of equivalence classes of matrices $\wt
W(\pi), \wt W_s(\pi)$
and $B(\pi)$ over $\Zpi$ that appear
in Theorem \ref{thm.cd}. These are mild variations of Witt-type
constructions,
motivated entirely by the geometric results of Section
\ref{sub.geometric}.

The group-ring $\Zpi$ has an involution defined by
$\overline{ng}=ng^{-1}$
for $n \in \BZ$ and $g \in \pi$. Let $A$ be a {\em Hermitian matrix} 
over $\Zpi$, ie, one that satisfies $\bar{A}^t=A$, where ${}^t$
denotes
the transpose. Two Hermitian matrices $A,B$ are {\em congruent} if there exists
a non-singular matrix $P$ such that $B=PA\bar P^t$.  We say that a
Hermitian matrix $A$ is {\em almost even} if
for every $g\in\pi$ with $g^2 =1$ but 
$g\not= 1$, the coefficient of $g$ in any diagonal entry of $A$ is even.
Note that if $A$ is nonsingular and almost even, so is $A^{-1}$.
Also, any matrix congruent to $A$ is almost even.
Given two Hermitian matrices $A, B$, we will say that they are
{\em stably congruent} if there exist {\em unidiagonal} matrices $S_i$
such 
that the block sums $A\oplus S_1$ and $B\oplus S_2$ are congruent. 
A unidiagonal matrix is a diagonal matrix all of
whose diagonal entries are $\pm 1$. Let $\widetilde W (\pi )$ denote the
set 
of {\em stable congruence classes of non-singular almost even matrices}. 

\begin{proposition}\lbl{prop.grp}
$\widetilde W (\pi )$ is an abelian group under block sum.
\end{proposition}

\begin{proof}
We need to show that $A\oplus (-A)$ is stably congruent to a unidiagonal matrix.
In fact 
$A\oplus (-A)$ is
congruent to $\left(\smallmatrix 0&A\\A&A\endsmallmatrix\right)$, which
is 
congruent to
$\left(\smallmatrix 0&I\\I&A\i\endsmallmatrix\right)$ and which, in
turn, 
since $A\i$ is
almost even, is congruent to 
$\left(\smallmatrix 0&I\\I&D\endsmallmatrix\right)$, where $D$
is a diagonal matrix all of whose entries are $0$ or $1$. But this is 
congruent to some
unidiagonal matrix.
\end{proof}

\begin{remark}\lbl{rem.meta}
The proof shows that any {\em metabolic}  non-singular almost 
even matrix is trivial in $\widetilde W (\pi )$. A metabolic
matrix is one which is
congruent to a matrix of the form $\left(\smallmatrix
0&I\\I&X\endsmallmatrix\right)$, for
some $X$.
\end{remark}

\begin{remark}\lbl{rem.simp}
As a variation on this we can define $A, B$ to be {\em simple stably 
congruent} if $A\oplus S_1
=P(B\oplus S_2 )\bar P^t$, for some non-singular {\em elementary} matrix
$P$.
An 
elementary matrix is a product of matrices, each of which differs from
the identity matrix in one
of the two following ways: (i) there is a single non-zero off diagonal
entry, or (ii) one of the
diagonal entries is replaced by $\pm g$, for some $g\in\pi$.
If we then define
$\wt W_s (\pi )$ to be the set of {\em simple stable congruence classes
of 
elementary non-singular almost even matrices}, the same proof shows that 
$\wt W_s (\pi )$ is a group. There is an obvious homomorphism 
$\wt W_s (\pi )\to \wt W(\pi )$.
\end{remark}

\medskip
Let $B (\pi )$ denote the set of {\em simple stable congruence classes
of 
almost even non-singular Hermitian matrices}. This is a semigroup under 
block sum but is not a group since the proof of
Proposition \ref{prop.grp}, showing $-A$ is an inverse for $A$ under
simple 
stable congruence, 
 only works if $A$ is elementary---see Remark \ref{rem.simp}. There
is an 
obvious inclusion $\wt W_s (\pi )\sub B (\pi )$ and epimorphism 
$B (\pi )\to\wt W(\pi )$ whose composition
$\wt W_s (\pi )\to B (\pi )\to \wt W(\pi )$ agrees with the map
of Remark \ref{rem.simp}.

\subsection{Surgery and a link description of $\H(N)$}
\lbl{sub.geometric}

It is well-known that the set of closed 3--manifolds can be identified
with the set of framed links in $S^3$ modulo an explicit (Kirby) 
equivalence relation discussed below. The first goal of this section is to 
give a similar link description of the set $\H(N)$. The surgery obstruction 
maps $\Phi, \Phi_0$ and $\Phi_0^s$ will then be obtained by considering 
linking matrices of appropriate classes of links.

\begin{lemma}
\lbl{lem.hndl}
If $f\co M\to N$ is a degree $1$ map, then one can adjoin handles of index
$2$ 
to $N$ to obtain a compact $4$--manifold $V$ with $\bd V-N=M$, and extend
$f$ 
to a map $F\co V\to N$ so that $F\vert N=\text{identity}$.
\end{lemma}

\begin{proof}
$f$ represents an element in $\Omega_3 (N)\iso H_3 (N)$ and so its
bordism class is determined by the degree of $f$. It follows that 
$f$ is bordant to the identity map of $N$. This gives us the manifold 
$V$ and map $F$ without the desired handlebody structure. We must 
eliminate the handles of index not equal to $2$.

First of all we can eliminate handles of index $0$ and $4$ in the usual
way. Now a handle of index $1$ represents a boundary connect sum with
$S^1\times D^3$. We can replace this with $D^2\times S^2$, thereby 
changing $V$. The only problem is how to replace the map $F$ on this 
altered piece. Since $F$ is a retract of $N$, $F\vert S^1\times D^3$ 
represents an element of
$F_{\ast}(\pi_1 (N))$. By sliding one foot of this $1$--handle around a 
representative of this element in $N\sub V$ we can arrange that $F\vert 
S^1\times D^3$ is null-homotopic. Thus, after replacing this $1$--handle 
with a $2$--handle, we can also replace $F$. To get rid of the $3$--handles 
we regard them as $1$--handles on $M$ and apply the same argument since 
$F_{\ast}\vert\pi_1 (M)$ is onto.
\end{proof}

Consider a manifold $V$ as in Lemma \ref{lem.hndl} and a choice of 
$2$--handles. These are attached along a framed link $L\sub N$. Since
$F\vert N=\text{id}$, the components of $L$ are null-homotopic. 
Conversely, given a framed null-homotopic link $L$ in $N$ we can
construct $V$ and then extend the identity map of $N$ over $V$ to obtain a
degree $1$ map $f\co N_L\to N$, where $N_L$ denotes the result
of surgery on $L$. The only indeterminacy in this construction is the 
choice of the extension $F$. There is no indeterminacy if we 
assume that $N$ is {\em prime}, ie $\pi_2 (N)=0$. In case $N$ is not
prime we proceed as follows. Consider the set $\Deg (N)$ of all {\em
diffeomorphism classes of 
 degree} $1$ maps $f\co M\to N$.   
We introduce an {\em equivalence relation} $\sim$ in $\Deg(N)$ 
generated by the 
following modification of a map $f\co M\to N$.  Let $K\sub M$ be a framed knot 
representing an element in $\ker f_{\ast}\co \pi_1 (M)\to \pi_1 (N)=\pi$ 
and let  $\phi \co (S^2,\ast )\to (N ,x_0 )$ be any map. The framing of $K$ 
identifies a 
neighborhood $U$ of $K$ with $S^1\times D^2$ and we can assume that 
$f(U)=x_0$. Now define $f'\co M\to N$ by
$f'|\overline{N -U}=f|\overline{N -U}$ and $f'|U$ is given by the composition
$$
S^1\times D^2\overset{p_2}\longrightarrow D^2\overset{\rho}\longrightarrow
S^2\overset{\phi}\longrightarrow N
$$
where $p_2$ is projection on the second factor and $\rho$ is the 
identification map $(D^2 ,S^1)\to (S^2 ,\ast )$.

The above discussion gives a well-defined onto map 
$\L(N)\to \ti\Deg(N)$
of the set $\L(N)$ of {\em framed, nullhomotopic links}\footnote{Not 
to be confused with the notion of nullhomotopic links in the sense of Milnor.}
in $N$
to the set $\ti\Deg(N)$ of $\sim$--equivalence classes of $\Deg(N)$. This map 
is not one-to-one, since surgery on
nonisotopic links may correspond to the same element of $\ti\Deg(N)$.
Recall that two framed links $L, L' \subset N$ are {\em Kirby equivalent}
(which we denote by $L \simk L'$) iff $N_L$ is diffeo to $N_L'$.
Fenn--Rourke \cite[Theorem 8]{FR} showed that Kirby equivalence is generated
by three moves (and their inverses) on framed links:
\begin{itemize}
\item[(i)]
Adding a trivial knot with a $\pm 1$--framing in a ball disjoint
from the rest of the link.
\item[(ii)] 
Replacing a component with a connected sum of that component with a 
push-out of another component, suitably framed. 
\item[(iii)]
Adding a knot with arbitrary framing together with a meridian of it with
$0$--framing.
\end{itemize}
In a less advertised part of their paper, Fenn--Rourke considered the 
{\em equivalence relation} $\simkp$ on framed links generated by moves 
(i) and (ii) alone, see \cite[Theorem 6]{FR}, which is 
related to surgery on maps rather than surgery on spaces. 
It is easy to see that
$\simkp$ preserves the class of nullhomotopic links in $N$
and respects the map $\L(N)\to\ti\Deg(N)$.
A converse is given by consequence of a not-so-well known theorem of 
Fenn--Rourke \cite{FR}. 

\begin{proposition}\lbl{prop.fr1}
The map $\L(N)/(\simkp) \to \ti\Deg(N)$
is one-to-one and onto.
\end{proposition}

\begin{proof}
This follows from \cite[Theorem 6]{FR}. The isomorphism $i$ in  diagram
$(\Delta )$ of \cite{FR} is determined, since the maps $\pi_1
(M)\to\pi_1 (W({\bf L}_i))$ are isomorphisms. The commutativity of
$(\Delta )$ corresponds to the diffeomorphism equivalence of the pairs
$(N,f)$. 
\end{proof}

Given a framed link $L \in \L(N)$, we can define a {\em linking matrix}. For
each component $L_i$ of $L$ choose a lift $\ti L_i$ in 
$\ti N$, the universal cover of $N$.
Then we have linking numbers $\lk(g\ti L_i ,h \ti L_j )$ for any 
$g,h\in\pi$---if $i=j$ and $g=h$, we need to push off along some vector 
field in the given framing. Now the full linking element $\l (L_i ,L_j)
\in\Zpi$ is defined as $\sum_{g\in\pi}\lk(g\ti L_i, \ti L_j )g\i$.  
We will associate to $f$ the matrix $A=(\l (L_i ,L_j ))$.
There is a mild and manageable indeterminacy in the choice of lifts of $L$. In
particular, any change of lifts will change $A$ by a simple congruence.

\begin{proposition}
\lbl{prop.link}
The linking matrix $A$ is almost even. If $A$ is any almost even
Hermitian
matrix over $\Zpi$ then there is a framed link in $N$ with
null-homotopic 
components whose linking matrix is $A$.
\end{proposition}

\proof
Let $S$ be a surface in $\ti N$ bounded by $\ti L_i$. If $g\in\pi$ with 
$g^2 =1$ but $g\not=1$, then the coefficient of $g$ in $\l (L_i ,L_i )$
is 
the intersection number $S\cdot g\ti L_i$.  Consider the intersection 
$X=S\cap gS$. This is a collection of loops and arcs in $S$.
$\l (L_i ,L_i )$ can be computed by counting up (with sign) the number
of 
points of $X\cap\bd S$. Now $X$ is invariant under the action of $g$ and
any 
arc of $X$ with exactly one point on $\bd S$ is sent by $g$ to another
such 
arc in $X$. Since the action of $g$ is free this means there are an even 
number of such arcs.

To prove the realizability start out with any null-homotopic link in
$N$. 
Then the desired matrix $A$ can be obtained from the linking matrix of
this 
link by a sequence of two operations.
 \begin{enumerate}
 \item For any $g\in\pi$ and $i,j$ add $\pm g$ to $a_{ij}$ and
$\pm g\i$ 
to $a_{ji}$.
 \item For any $i$ add $\pm 1$ to $a_{ii}$.
 \end{enumerate}
 These operations can be realized by changing the link as follows.
 \begin{enumerate}
 \item Replace $L_i$ by a connected sum of $L_i$ with a small loop
linking 
$L_j$. The arc used
to take the connected sum is determined by $g$.
\item Change the framing of $L_i$ by introducing a single twist. \qed
\end{enumerate}

\begin{remark}
\lbl{rem.wall}
Another interpretation of the linking matrix $A$ is as a representative of 
the intersection 
pairing on $\ker\{F_{\ast}\co H_2 (\ti V)\to H_2 (\ti N)\}$. Note that 
$\ker F_{\ast}\iso H_2 (\ti V,\ti N)$ which is a free $\Zpi$--module with 
basis determined by the $2$--handles. $A$ is non-singular if and only if
$f$ 
is a $\Zpi$--homology equivalence and is, in addition, elementary if $f$
is a simple $\Zpi$--homology equivalence, see Wall \cite{Wa}. If we now narrow
the definition of the relation $\sim$ on $\H(N)$ by restricting our
modifications 
to knots $K$ satisfying $\l (K,K)=\pm 1$, in order to stay within the class of 
$\Zpi$--homology equivalences, then Proposition \ref{prop.fr1} implies 

\end{remark}

\begin{proposition}
\lbl{prop.fr2}
There is a one-to-one correspondence
$$
\Lns(N)/(\simkp)\to \ti\H(N)
$$
between the set of $\sim$--equival\-ence classes of $\H(N)$ and
the set of $\simkp$--equival\-ence classes of $\Lns(N)$ of framed nullhomotopic
links with a nonsingular linking
matrix.
\end{proposition}

Thus, we obtain a well-defined map $\Phi\co \H(N)\twoheadrightarrow
B(\pi)$ (which is onto by Proposition \ref{prop.link}) given by a composition 
$$\H(N)\to\ti\H(N)\cong \Lns(N)/(\simkp)
\to B(\pi)$$
which assigns to an element $f\co M\to N$ of $\H(N)$ represented by surgery
on a framed nullhomotopic link, the linking matrix of that link.

\begin{proposition}
\lbl{prop.stable}
If $f\co M\to N$ is a $\Zpi$--homology equivalence then the stable
congruence 
class of $A$ depends only on the $\Zpi$--homology bordism class of $f$.
\end{proposition}
\begin{proof}
Suppose that $f'\co M'\to N$ is bordant to the identity on $N$ by $F\co V'\to
N$, 
where $V'$ consists of $2$--handles adjoined to $N$. Let $A'$ be an 
associated matrix. Suppose $f'$ is bordant to $f$ by a $\Zpi$--homology 
bordism $G\co W\to N$. By pasting $V,V',W$ together we
create a bordism $\hat G\co X\to N$ from the identity map on $N$ to itself.
The 
intersection pairing on $\ker\hat G_{\ast}\co H_2 (\ti X)\to H_2 (\ti N)$
is 
represented by $A\oplus(-A')$. Now suppose $\hat G$ is bordant, rel
boundary, 
to the projection $I\times N\to N$. Then the standard argument shows
that 
the intersection pairing on $\ker\hat G_{\ast}$ is metabolic. Thus by 
Remark \ref{rem.meta} the proposition is proved. The obstruction to this
bordism is an element of the bordism group $\Omega_4 (N)\iso H_4 (N
)\oplus\Omega_4$. Now $H_4 (N)=0$ and $\Omega_4$ is generated by $\BC
P^2$ and 
so this bordism will exist after we connect sum, say, $V$ with a number
of 
copies of $\pm \BC P^2$. But this can be achieved by adding to the
framed 
link defining the handlebody decomposition of $V$ a number of trivial 
components with $\pm 1$--framing. The effect of this is to block sum $A$
with a
unidiagonal matrix.
\end{proof}

Thus, we have a well-defined map $\Phi_0 \co \H_0 (N)\to\wt W (\pi )$ which 
is onto by Proposition \ref{prop.link}. We could also construct an
analogous 
map $\Phi_0^s \co \H_0^s (N)\to\wt W_s (\pi )$. 
The commutativity of the diagram of Theorem \ref{thm.cd} is obvious.

\smallskip
\begin{example}
Consider $N =S^2\times S^1 ,\pi =\Z$. For example if $M$ is obtained by
$0$--surgery on a knot $K$  in a homology $3$--sphere, then there is an obvious
degree $1$ map $f\co M\to N$, which is a $\Z\pi$--homology equivalence if and only
if
$K$ has Alexander polynomial $1$.

Now suppose $f,g\co M\to N$ are $\Z\pi$--homology equivalences. If $f_* =g_* \co H_1
(M)$ $\to H_1 (N)$ then it follows from the Hopf classification theorem that
$f\vert M-\text{point}\simeq g|N-\text{point}$ (homotopic). Thus $f\simeq g\# h$
for some $h\co S^3\to S^2\sub S^2\times S^1$.  It follows, using the geometric
definition of the Hopf invariant,  that $f\sim g$. If $f_*\not= g_* $, then $f_*
=g'_*$, where $g'=r\circ g$ and $r$ is the
self-diffeomorphism of $S^2\times S^1$ obtained by reflecting both factors. 

This discussion shows that, for any $M$ which is $\Z\pi$--homology equivalent to
$N=S^2\times S^1$, there is either one or two equivalence classes of
$\Z\pi$--homology equivalences $M\to N$ depending on whether or not there is an
orientation-preserving self-diffeomorphism of $M$ which induces $-1$ on
$H_1(M)$.
\end{example}

\smallskip
\subsection{Comparison to surgery theory}
\lbl{sub.compares}
 
We explain here whyc $\Phi$ can be thought of, in a rough sense, as
encapsulating the {\em surgery-theoretic } invariants of $\H (N)$. This is not
meant to be a mathematically precise statement but more of a philosophical
statement. 

The surgery exact sequence of Browder--Novikov--Sullivan--Wall extends to
lower dimensions in the topological category according to Freedman--Quinn
\cite{FQ}. 
If $N$ is a closed oriented $3$--manifold and $\pi =\pi_1 (N)$ is 
{\em good} in the sense of Freedman--Quinn (which, admittedly, may be rare) then
we have an exact sequence:
$$
[\Sigma N : G/\text{Top}]\to L_0^h (\pi )\to\H_0^{\text{top}}(N)
\to [N :G/\text{Top}
]\to L_3^h (\pi )
$$
where $\H_0^{\text{top}}(N)$ is the topological version of $\H_0 (N)$
and $L_i^h (\pi)$ are the {\em Wall surgery obstruction groups} \cite{Wa}. 
It is known
that $[N :G/\text{Top}]$ can be identified with the set of topological
bordism classes
of degree $1$ maps $M\to N$ where all maps are equipped with a morphism
of the tangent bundles. We have isomorphisms 
$[N :G/\text{Top} ]\iso H^2 (N ;\Z/2 )$ and $[\Sigma N :G/\text{Top} ]
\iso H^1 (N ;\Z/2 )\oplus H^3 (N)$, where the $H^3 (N)\iso\Z$ summand is 
mapped isomorphically to $L_0^h (1) \sub L_0^h(\pi )$---see the survey
article of Kirby--Taylor \cite{KT}. Thus we get an
exact sequence
\begin{equation}\lbl{eq.surh}
H^1 (N ;\Z/2 )\to L_0^h (\pi )/L_0^h
(1)\overset{\tau_h}\longrightarrow\H_0^{\text{top}}(N)\to H^2 (N ;\Z/2 )
\to L_3^h
(\pi )
\end{equation}
Similarly we have an exact sequence
\begin{equation}\lbl{eq.surs}
H^1 (N ;\Z/2 )\to L_0^s (\pi )/L_0^s(1)\overset{\tau_s}\longrightarrow
\H_0^{s,\text{top}}(N)\to H^2 (N ;\Z/2 )\to L_3^s(\pi )
\end{equation}
for the analogous classification of simple homology equivalence.

There are obvious maps $L_0^h (\pi )\to\wt W(\pi ), L_0^s (\pi )\to\wt
W_s 
(\pi )$. From the definition of stable congruence these maps induce maps 
$L_0^h (\pi )/L_0^h (1)\to\wt W(\pi )$ and $L_0^s (\pi )/L_0^s (1)\to
\wt W_s(\pi )$. According to \cite[Prop. 8.2]{R} this is an isomorphism
modulo $8$--torsion. It is not generally an isomorphism--see \cite[Theorem 10.4]{R}
for an example with $\pi =\Z\times\Z$. 

If $\pi$ were good, and we were able to ignore the difference between smooth and
topological equivalence, we could think of the maps $\Phi_0$ and $\Phi_0^s$ as 
approximations to left inverses of the maps $\tau_h$ and $\tau_s$ from the 
sequences \eqref{eq.surh} and \eqref{eq.surs}. Thus  we can roughly think of
$\ker\Phi$ 
as those manifolds $\Zpi$--homology equivalent to $N$ which
are undetected by surgery theory. In our theory of finite-type 
invariants, the invariants of degree $0$ will be those which can be
recovered 
from surgery theory, in this sense,  and conventional algebraic topology. Those
of
positive 
degree can detect differences invisible to surgery theory and
conventional
algebraic topology.

\section{Finite type invariants}
\lbl{sec.fti}

\subsection{$\K(N)$ and finite type invariants}
\lbl{sub.KtoM}

In this section we study the kernel $\K(N)$ of the map $\Phi$, which
leads rather naturally to a distinguished class of links in $N$ and to
a notion of \fti s. 

Suppose $f\co M\to N$ represents an element of 
$\K (N)$. Then, by Lemma \ref{lem.hndl}, there is a framed link $L\sub
N$ with null-homotopic components and linking matrix 
simply stably congruent to a unidiagonal matrix, such that if $V$ is obtained
from $N$ 
by adding handles along $L$, then $M=\bd V-N$ and the identity map of $N$ 
extends to a map $F\co V\to N$ so that $F|M=f$ . 
Since the moves that define simple stable congruence can be realized 
by either handle slides, choosing a different lift of $L$
or adding a trivial $\pm 1$--framed component to $L$ we  can, in fact,
assume that the linking matrix of $L$ is unidiagonal. We will call a 
(framed) link in $N$ whose components are null-homotopic
and linking matrix is unidiagonal {\em $\pi$--algebraically split}
($\pi$--AS
in short).  Conversely, given a $\pi$--AS link $L$ in $N$ we can
construct $V$ and then extend the identity map of $N$ over $V$ to obtain
an 
element of $\K(N)$. Let $\ti\K(N)$ denote the set $\K(N)/(\sim)$ for 
the equivalence relation $\sim$ of Section \ref{sub.geometric} and 
let $\Las(N)$ denote the set of framed $\pi$--AS links in $N$. It follows from 
Proposition \ref{prop.fr1} that $\ti\K(N)$ is in
one--one correspondence with the set $\Las(N)/(\simkp)$.
Note that a handle-slide will usually not preserve the property of being
$\pi$--AS, so our moves will be sequences of Kirby moves. 

We now imitate the usual approach that defines a notion of \fti s on a 
set of objects equipped with a move. Let $\F (N)$ denote the free abelian
group on the set $\ti\K(N)$. For a $\pi$--AS link $L\sub N$ we define 
$[N ,L]=\sum_{L'\sub L}(-1)^{|L'|}N_{L'}\in\F (N)$, where $N_L$ denotes
the result of surgery of $N$ along $L$ and $|L|$ denotes the number of
components of $L$. Note that the accompanying map $f'\co N_{L'}\to N$ is
uniquely determined by $L'$, modulo $\sim$, so we will
 suppress $f'$ from the notation. There is a decreasing filtration $\F
(N)=\F_0(N)\supset
\F_1(N) \supset \F_2(N)\dots$ on $\F(N)$, where $\F_n (N)$ 
denote the subgroup of $\F (N)$ generated by all $[N ,L]$ with $|L|\ge
n$. 
We call a function $\l\co  \ti\K(N)\to A$ with values in an abelian group $A$
a {\em \fti \  of type} $n$ iff $\l(\F_{n+1}(N))=0$.

For a decreasing filtration $\F$ (such as $\Fas {}$ or $\FY {}$ below)
we  let $\G$ denote the graded quotients defined by
$\G_n=\F_n/\F_{n+1}$. 

\subsection{Functoriality}
\lbl{sub.functoriality}

Since surgery theory is functorial, we might expect this to also be true of 
its deformation given by \fti s.

Suppose $f\co N'\to N$ is a $\Zpi$--homology equivalence between 
two closed oriented $3$--manifolds, where $\pi =\pi_1 (N)$. Then 
$f\in\H (N)$. Consider the induced function $f_{\ast}\co 
\H (N' )\to\H (N)$, where $f_{\ast} (g)=f\circ g$, and the further induced
function  $f_{\ast}\co  \ti\H (N' )\to\ti\H (N)$.
We will need the following lemma.

\begin{lemma}\lbl{lem.null}
 Suppose $K,L$ are disjoint links in $N$ such that each component of $K$
is null-homotopic in
$N$. Then $K$ is isotopic in $N$ to a link $K'$ such that each component
of $K'$ is null-homotopic in $N-L$. 
\begin{enumerate}
\item If the components of $L$ are also
null-homotopic in $N$, then the linking elements 
$\l (K'_i ,L_j )=0$, for every pair of components of $K' ,L$.
\item If $K$ is algebraically split in $N$, then we can choose $K'$ so that it
is algebraically split in $N-L$.
\end{enumerate}
\end{lemma}

\begin{proof} 
Since each component $K_i$ of $K$ is null-homotopic in $N$, it is
homotopic in $N-L$ to a
product of meridians of $L$. Thus we can connect sum several meridians
of $L$ to $K_i$ to get a
new knot $K'_i$ which is null-homotopic in $N-L$ and is clearly isotopic
to $K_i$ in the complement of the other components of $K$.

To see that $\l (K'_i ,L_j )=0$, when  $L_j$ is null-homotopic, we only need
note that any lift $\ti
K'_i$ to the universal cover $\ti N$ of $N$ is null-homotopic in 
$\ti N-\ti L_j$, for any lift of $L_j$.

If $K$ is algebraically split then the linking elements $\l (K'_i ,K'_j
)\in\BZ\pi$, where $\pi =\pi_1 (N-L)$, differ from the entries of a unidiagonal
matrix by members of the two-sided ideal $I$ of $\BZ\pi$
generated by elements of the form $h-1$, where $h\in G=\ker\{\pi_1 (N-L)\to\pi_1
(N)\}$. Figure \ref{figlink} shows how to modify $K'$ to change $\l (K'_i ,K'_j
)$ by an element:
\begin{enumerate}
\item  $g_1 (h-1)g_2$ if $i\not= j$,
\item $g_1 (h-1)g_2 +g_2\i (h\i -1)g_1\i$ if $i=j$
\end{enumerate}
for $g_1 ,g_2\in\pi ,h\in G$, without changing
any other linking element $\l (K'_r ,K'_s )$ except when  $r=j,s=i$. 
\begin{figure}[ht!]
\psdraw{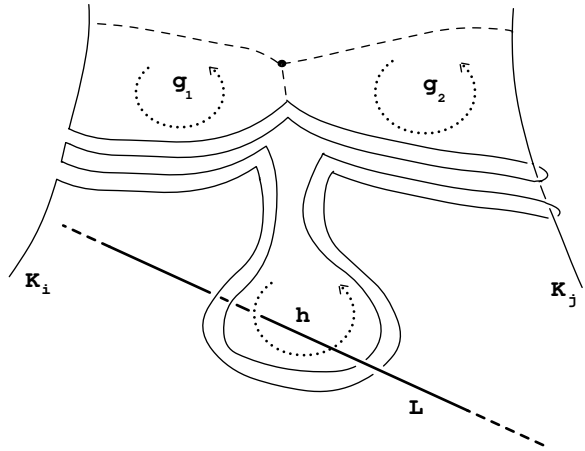}{3in}
\caption{Modification of a link}\lbl{figlink}
\end{figure}

The dotted
curves connecting $K'_i ,K'_j$ to the basepoint are those used to specify the
lifts---these are needed to define the linking elements. Note that we can
represent $h$ by the boundary of a disk in $N$ which is disjoint from $K'$ and
the arcs used in the modification, since $h$ is a product of meridians of $L$.
Thus the modified $K'$ is isotopic to $K$ in $N$. 

Since elements of the form (1) generate $I$, we only need show, by Proposition
\ref{prop.link}, that
self-conjugate almost even elements of $I$ are linear combinations of elements
of the form
(2), which we will call {\em norm-like}, to conclude that $K'$ can be chosen to
be algebraically split. 

Choose a subset $S\sub\pi_1 (N)$ so that, for every $g\in\pi_1 (N)$, exactly one
of $g,g\i$ belongs to $S$. For each $g\in S$ choose $\ti g\in\pi$ so that $\ti
g\to g$; choose $\ti 1=1$. Now suppose $\l$ is a self-conjugate even element of
$I$.  We can write $\l$ uniquely in the form
\begin{equation}\lbl{eq.diag}
 \l =\sum_{g^2 \not= 1}(\l_g\ti g +\ti g\i\bar\l_g )+\sum_{g^2 =1}\l_g\ti g 
 \end{equation}
 where $\l_g\in IG$, the augmentation ideal of $\BZ G$. Clearly the terms of the
first summation in equation \eqref{eq.diag} are norm-like, so we consider each
term $\l_g\ti g$ of the second summation.  Let us write $\l_g =\sum_i\e_i g_i$
and so
\begin{equation}\lbl{eq.ord2}
\l_g\ti g=\sum_i\e_i g_i\ti g
\end{equation}
where
$g_i$ are distinct elements of $G$, $\e_i\in\BZ$  and $\sum_i\e_i =0$.  Since
$\l$ is self-conjugate we have
$\l_g\ti g =\ti g\i\bar\l_g$ and so, for each $i$ there is some $j$ so that
$\e_i g_i\ti g=\e_j\ti g\i g_j\i$. If $i\not= j$, then replace $\e_j g_j\ti g$
in equation \eqref{eq.ord2} by $e_i\ti g\i g_i\i$. If $i=j$ then $g_i\ti g$ is
of order $2$ and so
$\e_i$
is even. In this case rewrite $\e_i g_i\ti g$ in equation \eqref{eq.ord2}, as
$2\e_i g_i\ti g$. Now equation \eqref{eq.ord2} will look like
\begin{equation}\lbl{eq.norm}
\l_g\ti g =\sum_i\e_i (g_i\ti g +\ti g\i g_i\i )
\end{equation}
where still $\sum_i\e_i =0$. If we now subtract $0=\sum_i\e_i (\ti g+\ti g\i  )$
from equation \eqref{eq.norm} we get 
$$
\l_g\ti g =\sum_i\e_i ((g_i -1)\ti g +\ti g\i (g_i\i -1))
$$
which is a sum of norm-like terms.
\end{proof}
\begin{proposition}\lbl{prop.phi}
We have:
$$\Phi (f_{\ast}(g))=\Phi (f)+f_{\ast}\Phi (g). $$
\end{proposition}
\begin{proof} 
Suppose $L\sub N$ is an algebraically split link determining $(N'=N_L , f)$ 
in $\H(N)$ and $K\sub N'$
determines $(M=N'_K ,g)$ in $\H (N')$. We can apply Lemma \ref{lem.null} to $K$
and the meridians $L'$
of $L$ in $N'$ to allow us to assume that the components of $K$ are
null-homotopic in
$N'-L'=N-L$. Now $K\cup L$ is algebraically split in $N$ and determines
$(M, f\circ g)$ in $\H
(N)$. Since, by Lemma \ref{lem.null}, $\l (K_i, L_j )=0$, the linking
matrix of $K\cup L$,
which  represents $\Phi (f\circ g)$ is the block sum of the linking
matrix of $L$, which
represents $\Phi (f)$ and the image under $f_* \co \pi_1 (N')\to\pi$ of the
linking matrix of $K$, which represents $f_* \Phi (g)$.
\end{proof}

\noindent
It follows from Proposition \ref{prop.phi} that $f_\ast
(\ti\K(N'))\sub\ti\K(N)$. 

\begin{proposition}
\lbl{thm.ihs}
$f_{\ast}(\F_n (N' ))\sub\F_n (N )$, for any $n$, and the induced maps
$f_{\ast}\co \G
(N')\to\G
(N)$ are epimorphisms.
\end{proposition}

\begin{proof}
 Suppose that $L\sub N$ is an algebraically split link which defines
$(N',f)$. Now let $K\sub N'=N_L$ be any algebraically split link---XSwe
can assume that $K$ is disjoint from the meridians of $L$ and so lies in
$N'-L'=N-L$. In fact, by Lemma \ref{lem.null} applied to $K$ and the
meridians of $L$,  we can assume that the components of $K$ are
null-homotopic in $N-L$ and that $K\cup L$ is algebraically split in
$N$. It is clear that $K\cup L$ defines the element $f_* (M,g)\in\K
(N)$.

Suppose that $K$ has $n$ components and so $[N',K]\in\F_n(N')$. Then we
have:
\begin{equation*}
f_*([N',K]) =[N_L,K]=\sum_{L'\sub L}(-1)^{|L'|}[N, K\cup L']
\end{equation*}
Thus we see that 
$f_* ([N', K])\in\F_n (N)$, which shows that $f_{\ast}(\F_n (N' ))
\sub\F_n (N )$, and that
$f_* ([N', K])\con [N, K]\mod\F_{n+1}(N)$, which shows that
$f_{*,n}\co \G_n(N')$ $\to\G_n (N')$ is onto.
\end{proof}

We also have natural filtration-preserving maps $\F (N)\to\F (N\# N'
)$, for any $N , N'$, defined by
$$
f\co M\to N\quad\rightsquigarrow\quad f\#\text{id}\co M\# N'\to
N\# N'.
$$
In particular $\F(N)$ is an $\F(S^3)$ module via a map $\F (S^3 )\to\F
(N)$, 
for any $N$.  

\section{The $\Fas {}$ and $\FY {}$ filtrations}
\lbl{sec.compare}

In this section we reformulate our theory of \fti s in terms of
\clover s which in particular allows us to deduce upper bounds for the
number of invariants in terms of $\pi$--decorated trivalent graphs.
Recall the notion of a {\em \ygraph \ } in $N$ from \cite{Gu,GGP},
the terminology of which we follow here. A {\em \ylink \ } is a disjoint 
union of \ygraph s, and a {\em \clover \ } is a mild generalization
of a \ylink . Given a \clover \ $G$ in $N$, let $N_G$ denote the result
of surgery on $G$. Throughout this paper, by a {\em \ylink} or 
{\em clover} in a manifold $M$ we will mean one with {\em nullhomotopic
leaves}. This condition, dictated by our surgery problem of Section 
\ref{sec.surgery}, matches perfectly the generalization of the results of
\cite{GGP} from the case of $N=S^3$ to the case of arbitrary $N$.

Recall that surgery on a \ygraph \ is equivalent to surgery
on a six component framed link which consists of the three edges and the 
three leaves of the \ygraph .
Since surgery on a \ygraph \ $G$ (or more generally, a clover)
is an example of surgery on a nullhomotopic link with non-singular linking
matrix, it follows that $N_G \in \ti\H(N)$.  An alternative geometric proof 
of this may be obtained from the fact that $G$ lifts to $\pi$ copies of 
\ygraph s $\ti G$ in $\ti N$ and that the $\pi$ covering $\ti {N_G}$ can be 
identified with $(\ti N)_{\ti{G}}$---since surgery on \ylink s in a 3--manifold 
does not change its homology, it follows that $(\ti N)_{\ti{G}}$ 
is $\Zpi$--homology cobordant to $\ti N$. Since the linking matrix of a 
\ygraph \ is a metabolic matrix, it follows by Remark \ref{rem.meta} that
$N_G \in \ti\K(N)$. 

Let $\ti\KY(N)$ denote the subset of
$\ti\K(N)$ that consists of all maps $N_G$ for clovers $G$ in $N$,
and let $\FY {}$ denote the free abelian group on $\ti\KY(N)$.
We define a decreasing filtration $\FY {}=\FY 0 \supset \FY 1
\supset \FY 2 \dots$ on the abelian group $\FY {}$
where $\FY {n}$ denotes the span of $[M,G]$ for clovers $G$ in $M$ of
degree
(ie, number of trivalent vertices) at least $n$ with nullhomotopic
leaves. 

We will show later that $\ti\K(N)=\ti\KY(N)$ and that for all 
integers $n$, we have $\FY {2n}=\Fas {3n}$ after tensoring with $\BZ[1/2]$.

\subsection{The $\A$--groups and the graded quotients $\GY {}$}
\lbl{sec.gradedq}

The discussion of \cite[Section 4.3]{GGP} implies that $\GY n 
\oZ$ is generated by $[N,G]$ for clovers of degree $n$ without leaves,
ie,
embedded trivalent graphs of degree $n$. Unlike the case of $N=S^3$,
however,
$[N,G]$ depends on the embedding. Consider two embedded trivalent graphs 
$G, G'$ in $N$ such that $G\sminus e=G'\sminus e'$ for edges $e,e'$ of
$G,G'$ which are homotopic, rel. boundary. Then, the Sliding Lemma (in
the 
form of \cite[Corollary 4.2]{GGP}) implies that $[N,G]=[N,G'] \in \GY
{}$.
The following lemma describes the induced equivalence relation on 
embedded trivalent graphs.

\begin{lemma}
\lbl{lem.hclasses}
For an abstract (not necessarily connected) graph $G$ and path 
connected space $X$ we have a 1--1 onto map
$$
[G,X]\cong \out(\pi_1(G), \pi_1(X))
$$
where $\out(\pi_1(G), \pi )=\prod_{\a}\out(\pi_1(G_{\a}), \pi_1(X))$, the
Cartesian product over the connected components of $G$, and $$\out (G_1 ,G_2
)=\Hom (G_1 ,G_2 )/(\mathrm{inner\ automorphisms\ of}\ G_2 ).$$
\end{lemma}

\begin{proof}
Pick a {\em maximal forest} $T$ for $G$, ie, a maximal tree for each 
connected
component of $G$. If $f \in [G,X]$, then we can assume that $f(T)=x_0$, 
a base point of $N$, ie, that it factors though a map $G/T\to X$. This
map is determined by the induced one on the level of $\pi_1$. A
different
choice of a maximal forest or a different choice of a base point of $N$
results in maps on $\pi_1$ that differ in each connected component of $G$
by independent inner automorphisms of $\pi_1(X)$.
\end{proof} 

Let $\A'(\pi)$ denote the abelian group generated
by pairs $(G,\a')$ for abstract (not necessarily connected)
vertex-oriented trivalent graphs $G$ together with
an $\a' \in \out(\pi_1(G),\pi)$,
modulo the $\AS$ and $\IHX$ relations. We call
$\a'$ a {\em $\pi$--decoration} of $G$.
For each pair $(G,\a')$, pick an arbitrary embedding of $G$ in $N$
so that the induced map on the fundamental group coincides with $\a'$.
Equip the embedding with an arbitrary framing, thus
resulting in a clover in $N$. \cite[Lemma 4.4, Corollary 4.5 and 
Theorem 4.11]{GGP} shows that this define
a map $\Psi_n\co \A'(\pi)\to \GY n$. The above discussion implies that:

\begin{theorem}
\lbl{thm.graded}
For every $n$, the map $\Psi_n\co   \A'_n(\pi)\oZ \to \GY n \oZ$
is onto and functorial with respect to $\Zpi$--homology equivalences.
\end{theorem}

\begin{remark}
Over $\BQ$, and for $\pi_1(N)=1$, it is known that the map $\Psi$ is an
isomorphism, due to the existence of sufficiently many invariants
first constructed by Le--Murakami--Ohtsuki \cite{LMO}. 
\end{remark}

We now give an alternative description (closely related to $\pi$--AS
links, see Section \ref{sec.ssurgequiv}) of the notion of $\pi$--decoration
of a graph. This description generalizes to a decoration of the edges
by elements of an arbitrary ring with involution, see Definition 
\ref{def.Agroups} below.

Consider pairs $(G,\a)$ of abstract, vertex-oriented, {\em edge-oriented}
trivalent graphs $G$, together with a map $\a\co \edge(G)\to \Zpi$
that colors each oriented edge of $G$ by an element of $\Zpi$. 
Let $\A(\pi)$ denote the abelian group generated by pairs
$(G,\a)$ modulo the relations shown in Figure \ref{relations}.

\begin{figure}[ht!]
\psdraw{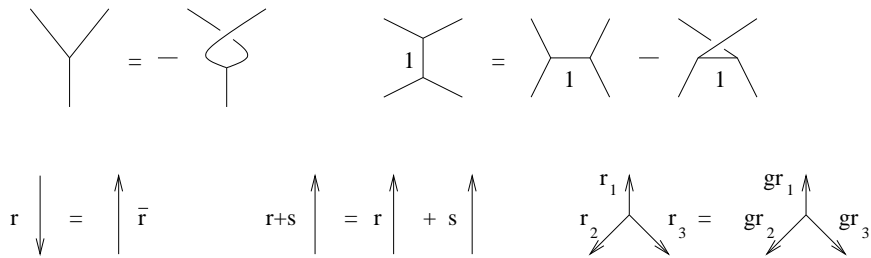}{4.5in} 
\caption{The $\AS$, $\IHX$, $R_1, R_2$ and $R_3$ relations. 
Here $\bar{g}=g^{-1}$ is the involution of $\Zpi$, $r,s,r_i \in \Zpi$
and $g \in \pi$.}\lbl{relations}
\end{figure}

\begin{lemma}
\lbl{lem.equal}
There is an isomorphism
$$
\A(\pi)\to\A'(\pi).
$$
\end{lemma}

\begin{proof}
It suffices to consider  a vertex-oriented, edge-oriented {\em connected} 
graph $G$. To a map $\a$, we will associate a map $\a'$ and vice versa.

Given a map $\a\co \edge(G)\to
\pi$, (which in view of relation $R_2$ we may assume that it is a decoration
of the edges of $G$ by elements of $\pi$)
we define a map $\a'\co \pi_1(G)\to\pi$ as follows.
For a closed path of oriented edges $e=(e_1,\dots,e_n)$, we set
$\a'(e)=\a(e_1)\dots\a(e_n)$. It is easy to see that this defines a
group homomorphism $\pi_1(G)\to\pi$, compatible with the relations $R_1$ 
and $R_3$. Conversely, given a map $\a'$, choose a maximal tree $T$ and define 
$\a(\edge(T))=\{1\}$. Since $\pi_1(G)$ can be identified with the
free group on  $\edge(G\sminus T)$, $\a'$ will then determine $\a$ on these
edges.
\end{proof}

The following concept of the $\A$--{\em groups}, motivated by Theorem 
\ref{thm.graded}, has several applications which will be presented in a later
publication.

\begin{definition}
\lbl{def.Agroups}
Given a ring $R$ with involution and a subgroup $U$ of its group of units, 
we define $\A(R,U)$ to be the graded abelian group generated by trivalent 
graphs (with a vertex and an edge orientation) whose edges are decorated by 
elements of $R$, modulo the relations of Figure \ref{relations}, with 
$r,s,r_i \in R$ and $g \in U$.
\end{definition}


\subsection{The equivalence of the $\Fas {}$ and $\FY {}$ filtrations}

In \cite[Sections 5.2--5.6]{GGP} it was shown that \fti s of 
\ihs s based on surgery on AS links coincide with those based on surgery
on clovers. In this section we will extend this to 3--manifolds, by using
the same idea as in \cite[Sections 5.2--5.6]{GGP}, together with Lemma 
\ref{lem.untie} and Proposition \ref{prop.untie}.

Recall that the proofs of \cite[Sections 5.2--5.6]{GGP} consist of three
types of arguments:
\begin{itemize}
\item
First, {\em untying} AS--links by \clover s and vice versa. 
\item
Second, counting arguments. 
\item
Third, an application of the topological
calculus of clovers to the study of $\GY {}$.
\end{itemize}
The second type of arguments works without change when we
replace
$S^3$ by $N$. So does the third type of argument, since, restricted
to clovers with nullnomotopic leaves, it uses the moves (i) and (ii)
of Kirby equivalence and the move (iii) for nullhomotopic knots, which
is a consequence of (i) and (ii) as shown in \cite[Theorem 2]{FR}.
The first type of argument requires some additional work, which
consists of Lemma \ref{lem.untie} and Proposition \ref{prop.untie}.

\begin{theorem}\lbl{thm.asY}
For all integers $n$ we have that $\Fas {3n}\oZ=\FY {2n}\oZ$.
\end{theorem}

Together with Theorem \ref{thm.graded}, it implies that:

\begin{theorem}
For all integers $n$ there is an onto map
$$
\A'_{2n}(\pi)\oZ\twoheadrightarrow \Gas {3n}\oZ
$$
functorial with respect to $\Zpi$--homology equivalences.
\end{theorem}

\subsection{Undoing \clover s by AS--links and vice versa}
\lbl{sub.untangle}

Let $G$ be a \ylink \ in $N$. Using the terminology of 
\cite[Section 5.3]{GGP}, we say that
a link $O$ in $N\sminus G$ {\em laces} $G$, if $O$ is
trivial, unimodular and each of the (pairwise disjoint)
discs bounding its components intersects $G$ in at most
two points, which belong to the leaves of $G$.
$G$ is {\em trivial}, if it consists of $n$ Y--graphs, 
standardly embedded in $n$ disjoint balls which lie
in an embedded ball in $N$.

\begin{lemma}
\lbl{lem.untie}\cite[Lemma 5.3]{GGP}
Let $T$ be a trivial $n$--component Y--link in $N$.
For any $n$--component Y--link $G$ in $N$, there
exists a unimodular link $O$ in $N$ which laces $T$, such that
$[N,G]=[N_O,T]$. Under surgery on $T' \subset T$,
$O$ gets transformed to a $\pi$--AS link in $N$.  
\end{lemma}

A \ylink \ $G$ in $M\sminus L$ {\em laces a link} $L$
if $L$ is $\pi$--AS, $G$ has nullhomotopic leaves
and every leaf $l$ of $G$ either bounds a disk
which intersects $L$ geometrically once, or the equivariant
linking number of $l$ and every component of $L$ vanishes.
We call $(G,L)$ a {\em lacing pair}.
A special case of a lacing pair $(G,L)$ for a trivial
\ylink \ $G$ was called a {\em Borromean surgery} in \cite{Ma} and a 
{\em $\Delta$--move} in \cite{MN}.

\begin{proposition}
\lbl{prop.untie}
Let $O$ be a trivial unimodular $n$--component link in $N$.
For any $n$--component $\pi$--AS--link $L$ in $N$, there 
exists a lacing pair $(G,O)$ such that $O$ is trivial unimodular,
under surgery on $G$ $(N,O)$ gets transformed into $(N,L)$
and under surgery on 
$O' \subset O$, $G$ gets transformed to a \ylink \ in $N$
with nullhomotopic leaves. 
\end{proposition}

\begin{proof}
Choose a base point $x_0$ of $N$ and a {\em basing} $\ga$ of $L$, ie,
a choice of disjoint paths $\{\g_i\}$ in $N\sminus L$ from $x_0$ to points 
$x_i\in L_i$, one for each component of $L$. Choose a framing of $L$ and 
a lift $\ti x_0 \in \ti N$ of $x_0$. Then, there is a unique
lift $\ti L \cup \ti \ga$ in $\ti N$ of $L \cup \ga$ that contains $\ti
x_0$ and a well-defined linking matrix $A$ of $L$.

Choose a regular homotopy $L^t$ from $L^0=L$ to $L^1=O$, which we can assume is
stationary on every $x_i$. We may assume that
$L^t$ is a link except for finitely many times $\{ t_s \}$ where
$L^{t_s}$
is an immersion with a single transverse double point. The linking matrix
$A_{-\e}$
and $A_{\e}$ of
$L^{t_s -\e}$ and $L^{t_s +\e}$ are related as follows: if the 
double point $p$ involves the components $L_i$ and $L_j$ of
$L^{t_s}$, with $i\le j$, 
construct a 
loop $g_p$ by starting at $x_0$, going along $\ga_i$ to $L_i$, then over
to
$L_j$ and back to $x_0$ via $\overline{\ga}_j$---see the figure below.

\psdraw{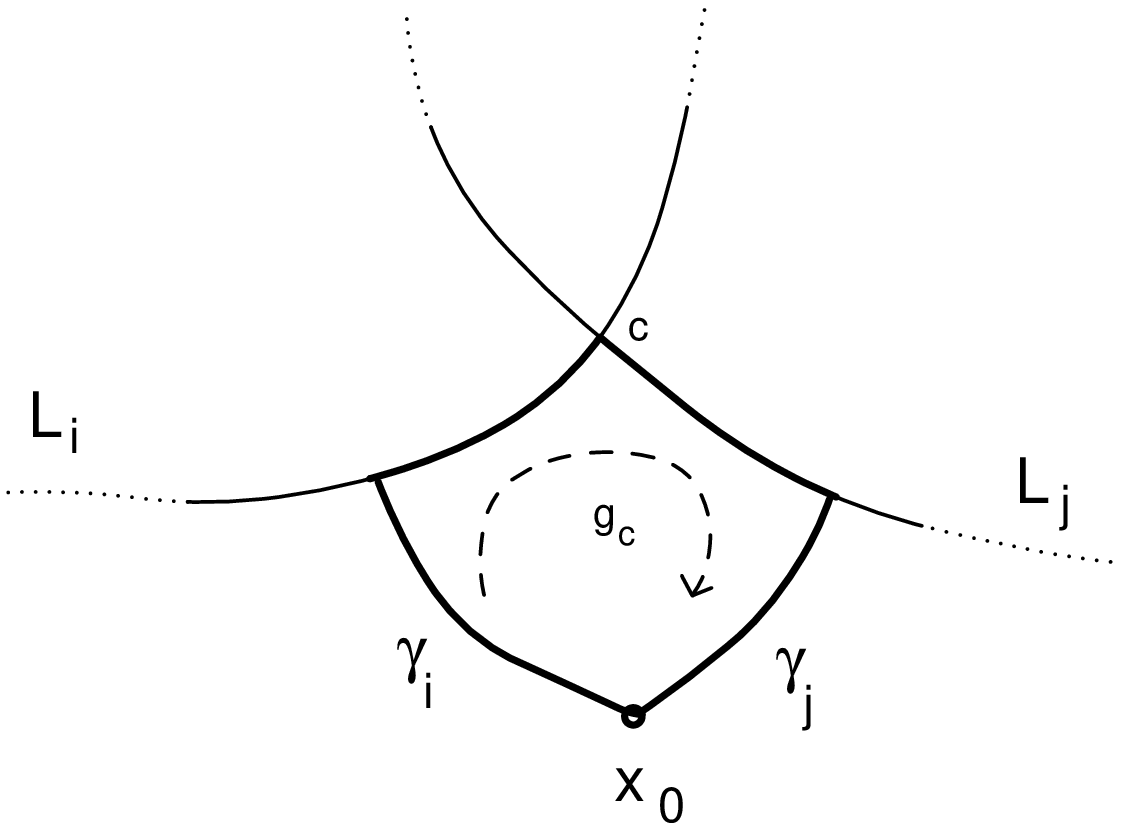}{2.5in}

It is easy to see that
$$A_{\e}=A_{-\e} + \e_p( g_p E_p+g_p\i E_p^t )$$
where 
$\e_p \in \{-1,+1\}$ is the local orientation sign of $p$ and $E_p$ is the 
matrix with all zeros except in the $(i,j)$ place where it equals 
$1$. Thus, if $A_L$ is the linking matrix of $L$ and $A_O$ is the linking 
matrix of an unlink with the same framing and number of components as $L$,
 we have 
$A-A_O=\sum_p\e_p( g_p E_p+g_p\i E_p^t )$, where the sum is over all
double points
of the homotopy. Since $L$ is a $\pi$--AS link, it follows that for every
pair of components $L_i$ and $L_j$ of $L$ there is a pairing of the double
points of $L_i$ and $L_j$ into classes $(p^+,p^-)$ such that
$g_{p^+}=g_{p^-}$ and $\e_{c^+} =-\e_{c^-}$. In other words,
one can undo $L$ by a sequence of {\em double crossing changes} shown
below in Figure \ref{fig.dblcr}, 
\begin{figure}[ht!]
\psdraw{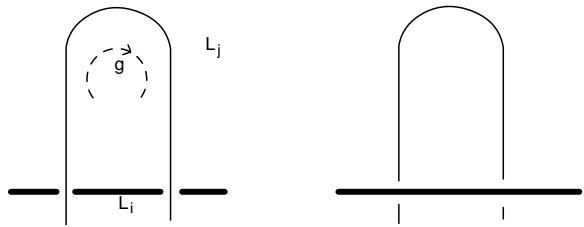}{3in}
\caption{A double crossing change}\lbl{fig.dblcr}
\end{figure}
where the loop $g\simeq g_{p^+}g_{p^-}\i=1 \in \pi$ is nullhomotopic.
These double crossing changes can be achieved by surgery on \ygraph s whose
leaves are nullhomotopic, 
see \cite{Ma,MN} and also  Lemma \ref{lem.bands} below. So far, each of
the \ygraph s have two leaves that bound a disk that intersects $L$ at most
once and a nullhomotopic leaf. Observe that every nullhomotopic leaf in $N$
bounds a disk with clasp intersections as shown below.
Using repeatedly Move $Y_4$ of \cite[Theorem 3.1]{GGP} (the so-called, move 
of {\em Cutting a Leaf}), as follows 
$$
\ppsdraw{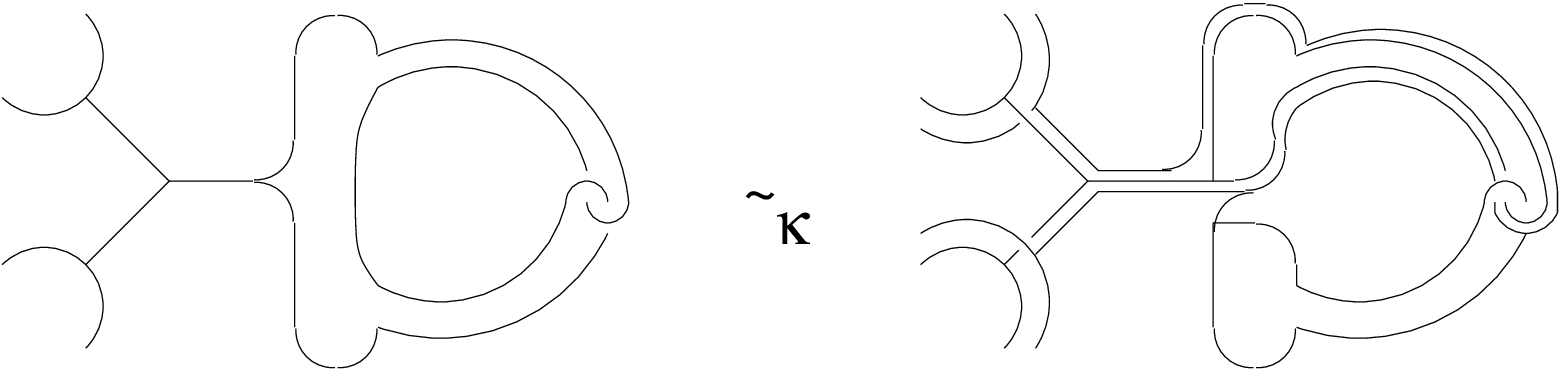}{4in}
$$
we may assume that every leaf of each \ygraph \ bounds a disk that
either intersects $L$ geometrically once, or none. In all cases, the
\ylink \ $G$ that consists of all these \ygraph s is lacing the unlink $O$. 
It is easy to verify that the rest of the statements of the proposition.
\end{proof}

The following lemma shows how to slide a band of an 
embedded surface through any collection of bands 
(or leaves of Y--graphs) by Y--surgery. 

\begin{lemma}\lbl{lem.bands}
The following framed links are Kirby equivalent:
$$\ppsdraw{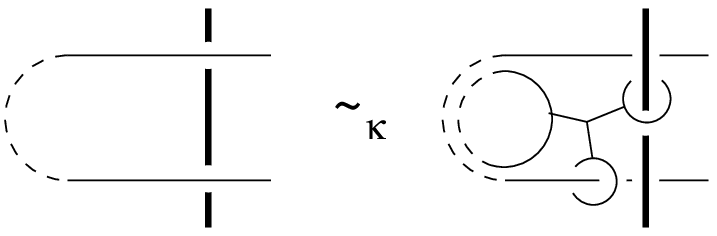}{2.5in}$$
In particular, a double crossing move can be obtained by surgery
on a \ygraph .
\end{lemma}

\begin{proof}
$$\ppsdraw{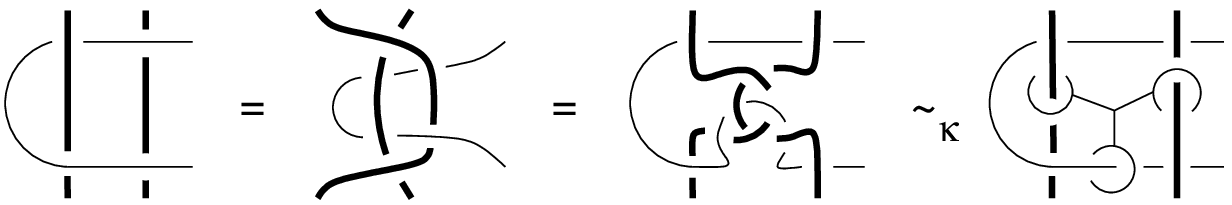}{4.5in}$$
\end{proof}

\begin{proposition}
\lbl{prop.K}
$\ti\K(N)=\ti\KY(N)$.
\end{proposition}

\begin{proof}
Since $\ti\KY(N)\subset \ti\K(N)$, we need only show the opposite
inclusion. Proposition \ref{prop.untie} implies that
for every $\pi$--AS link $L$ in $N$, there exists a trivial \ylink \
$T$ that ties a trivial unimodular link $O$ such that $(N,L)=(N_T,O)$.
Let $G$ denote the image of $T \subset N$ under surgery on $O$.
$G$ is a \ylink \ (with nullhomotopic leaves) and  
$N_L=N_{T \cup O}=N_{G}$. The result follows.
\end{proof}

\section{Surgery equivalence}
\lbl{sec.ssurgequiv}

In this section we discuss the notion of surgery
equivalence of $\pi$--AS links, motivated both by surgery theory
and by the theory of \fti s.  

Suppose that $L$ is an unframed $\pi$--AS link in $N$.
Now let $L$, with some unit framing, be expanded to a $\pi$--AS link $L\cup
L_{\text{triv}} \subset N$ 
for some  trivial, unit-framed link $L_{\text{triv}}$. Let $L'$
denote
the image of $L$ under the obvious isomorphism $N \cong
N_{L_{\text{triv}}}$.
{\em Surgery equivalence} is the relation on the set of unframed $\pi$--AS links
in $N$  generated by the
move that replaces  $L$ by $L'$ for some link $L_{\text{triv}}$ as above.

It was shown in \cite{Le} that, when $N =S^3$, surgery equivalence
classes of unframed 
AS--links are determined by the Milnor triple $\bar\mu$--invariants. We will
generalize the construction of these $\bar\mu$--invariants in an equivariant
manner to define surgery
equivalence invariants  of  $\pi$--AS links $L \subset N$. Choose a
lift 
$\ti L$ of $L$ in the $\pi$--cover $\ti N$. The components $\ti L_i$ of 
$\ti L$ bound oriented surfaces $V_i\sub\ti M$ and, since $L$ is 
algebraically split, we can assume that the interior of $V_i$ does not
intersect $p\i (L)$, where $p\co \ti N\to N$ is the projection. For any 
$1\le i,j,k\le q$ and
$g,h\in\pi$ we define $\bar\mu^{g,h}_{ijk}(L)$ to be the triple
intersection number of $V_i
, gV_j , hV_k$, when these are three different surfaces. This is
independent of the choice of $\{ V_i\}$. A change
in the choice of liftings $\{ \ti L_i\}$ produces the following change 
in $\bar\mu^{g,h}_{ijk}(L)$: for any $g_1, \ldots ,g_n\in\pi$
we have
$$\text{new }\{\bar\mu^{g,h}_{ijk}(L)\}=\text{old }\{\bar\mu^{\ g_i\i 
gg_j ,\  g_i\i
hg_k}_{ijk}(L)\}$$
Note that for the special case when $\pi$ is abelian, there is no 
indeterminacy in $\bar\mu^{g,h}_{iii}(L)$.

Also note the following:
\begin{enumerate}
\item $\bar\mu^{g,h}_{ijk}(L)=\bar\mu^{g\i h,g\i }_{jki}(L), \quad
\bar\mu^{h,g}_{ikj}(L)=-\bar\mu^{g,h}_{ijk}(L)$.
\item $\bar\mu^{g,h}_{ijj}(L),\ (i\ne j)$ is defined only if $g\ne h$.
\item $\bar\mu^{g,h}_{iii}(L)$ is defined only if $1,g,h$ are distinct.
\end{enumerate}

Now set
$$
\bar\mu_{ijk}(L)=\sum_{g,h\in\pi}\bar\mu^{g,h}_{ijk}(L)(g,h)\in
\Z [\pi\times\pi ]
$$
This is a finite sum.  In cases (2) and (3), the relations in (1) impose 
conditions on $\bar\mu_{ijj}(L)$. In case (2) it is skew-symmetric in
the two factors, but for case (3) the constraint is more complicated.

\begin{example}
If $\pi =\Z$, then $\bar\mu_{iii}(L)$ is determined by the 
$\{\bar\mu^{s,t}_{iii}(L)\ | s>t>0\}$. As an example of a knot $K$ with 
$\bar\mu^{s,t}_{111}=1$, choose a Borromean link
$L_1 ,L_2 ,L_3$ in a ball in $N$ and let $K$ be the connected sum 
$L_1\# L_2\# L_3$, where the tube connecting $L_1$ to $L_2$
winds 
$t$ times around the generator of $\pi$ and the tube connecting $L_2$ to 
$L_3$ winds $s-t$ times around. Note that $\bar\mu^{s',t'}_{111}(K)=0'$ 
unless $s=s'$ and $t=t'$. See the figure below for $s=t=1$.
\end{example}
$$
\ppsdraw{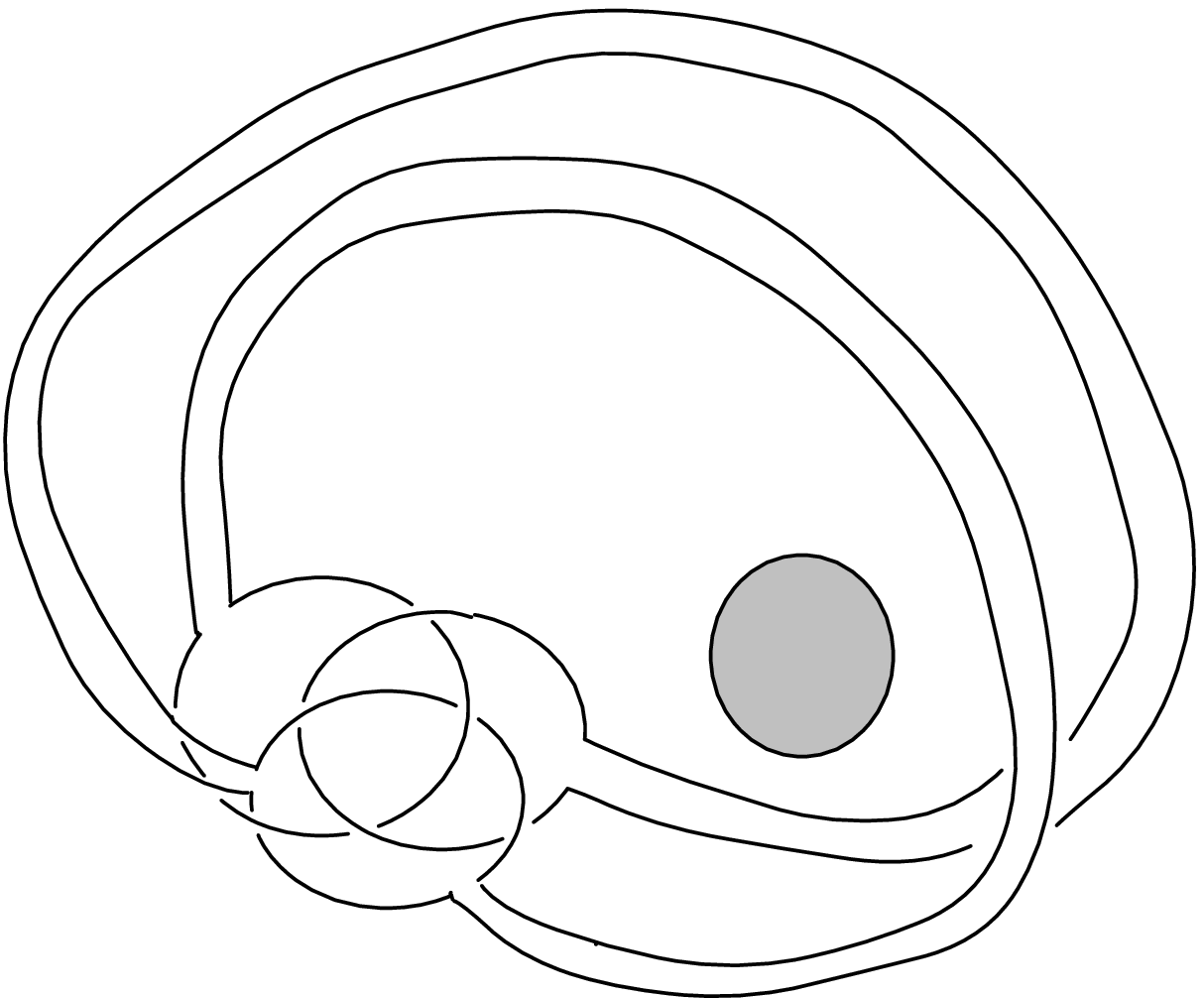}{2in}
$$
We will show that, just as in the case of unframed links in $S^3$, the surgery
equivalence class is determined by the $\mu_{ijk}^{g,h}$. First we need the 
following

\begin{lemma}
\lbl{lem.conc}
If two links $L$ and $L'$ differ as in the first two frames of this picture 
$$ 
\ppsdraw{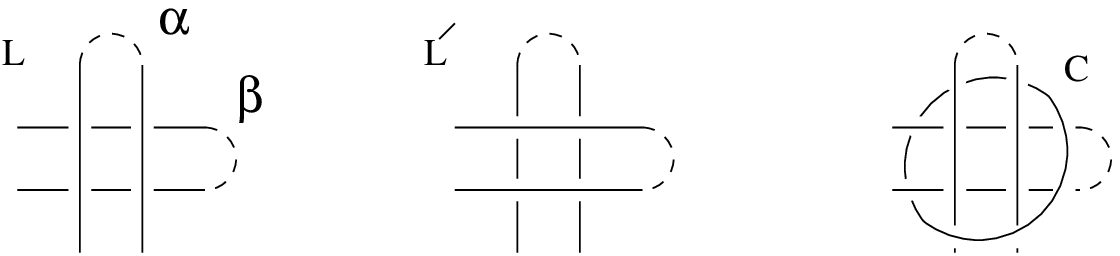}{3in} 
$$
or the first two frames of this picture
$$
\ppsdraw{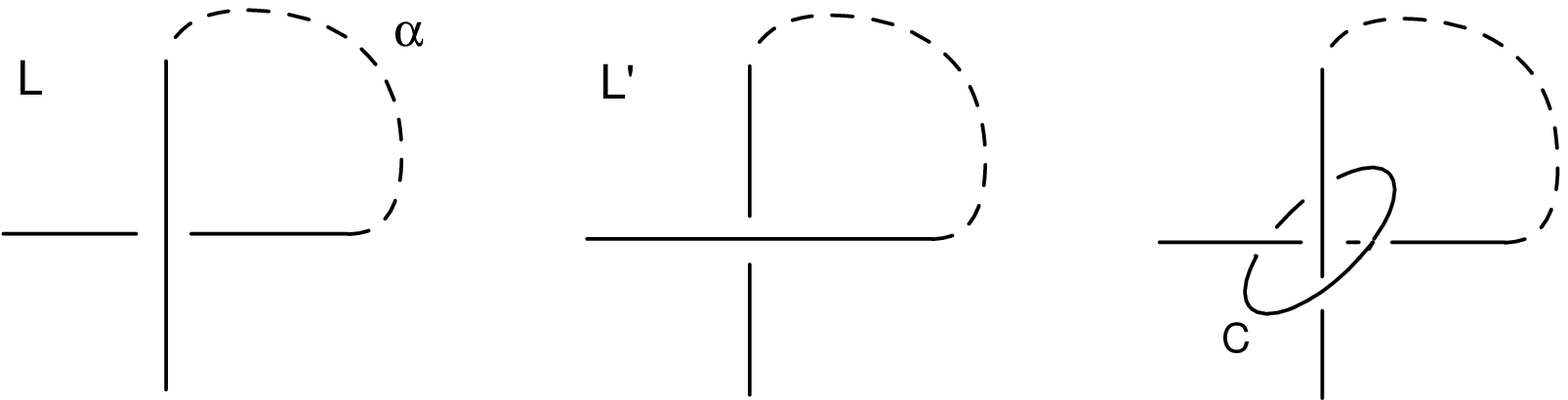}{3in}
$$
where $\a,\b$ are nullhomotopic paths, then they are surgery equivalent.
\end{lemma}

\begin{proof}
This follows by elementary properties of Kirby's calculus applied
to the unit-framed knots $C$ shown in the pictures above.
\end{proof}

\begin{theorem}\lbl{th.surg}  
Suppose $L$ and $L'$ are unframed $\pi$--AS links in $N$. Then $L$ is surgery 
equivalent to $L'$ if and only if $\mu_{ijk}^{g,h}(L)=\mu_{ijk}^{g,h}(L')$ 
for all $i,j,k,g,h$ for which they are defined.
\end{theorem}

\begin{proof}
First of all we prove the invariance of the $\mu_{ijk}^{g,h}(L)$ under
surgery equivalence. If $L_{\text{triv}}$ is chosen so that $L\cup
L_{\text{triv}}$ is $\pi$--AS, then we can, by tubing, arrange that the surfaces
$\{ V_i\}$ used to define $\mu_{ijk}^{g,h}(L)$ are disjoint from the lifts of
$L_{\text{triv}}$ and so pass unchanged into the $\pi$--covering of the surgered
link. In particular the intersections which define $\mu_{ijk}^{g,h}(L)$ are
unchanged. 

Now suppose that $L$ and $L'$ are two $\pi$--AS links such that
$\mu_{ijk}^{g,h}(L)=\mu_{ijk}^{g,h}(L')$. By Proposition \ref{prop.untie} we
know that $L$ can be transformed into $L'$ by surgery on a set
of \ylink s whose leaves are meridians of $L$. Since surgery on such a 
\ylink \ is the same as a sequence of disjoint {\em $\Delta$--moves} in the 
terminology of \cite{MN}---see Figure \ref{fig.delta}---it is easy to see the 
effect of such a surgery on the $\{\mu_{ijk}^{g,h}\}$.  
\begin{figure}[ht!]
\psdraw{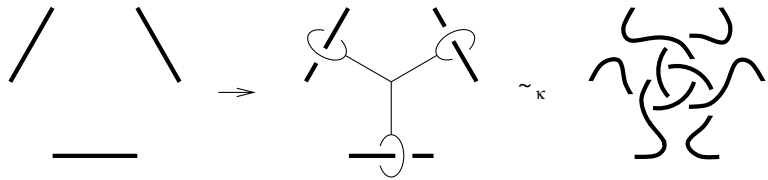}{4in} 
\caption{A $\Delta$--move}\lbl{fig.delta}
\end{figure}
Suppose a lift of the \ylink \ $G$ into $\ti N$ has its
meridians on three  components $\ti L_i , g\ti L_j ,h\ti
L_k$, where $i\le j\le k$. If any two of these components are the same then
there is no change in any of the $\mu_{i'j'k'}^{g',h'}$. If the three components
are distinct  then
$\mu_{ijk}^{g,h}$ is changed by $\pm 1$ and every other $\mu_{i'j'k'}^{g',h'}$,
where $i'\le j'\le k'$, is unchanged. Thus our assumption about $L, L'$ says
that the transformation from $L$ to $L'$ is accomplished by a sequence of
surgeries of two types:
\begin{itemize}
\item[(i)] surgery on pairs of
\ylink s $\{ G_i ,G'_i\}$, where $G_i$ and $G'_i$ can be lifted to \ylink s in
$\ti N$ with oppositely
oriented
trivalent vertices and which have leaves on the same three distinct components, 
and 
\item[(ii)] surgeries on individual \ylink s $G_j$ with at least two leaves
on the same component. 
\end{itemize}
In case (ii) it is easy to see that surgery on
$G_j$ does not change the surgery equivalence class since we can undo the
Borromean part of the \ylink\  by crossing changes, using the second part of
Lemma
\ref{lem.conc}, on the two rings attached to the same component of $L$. 

Thus it
remains to show that the effect of surgery on  a 
pair of \ylink s $G, G'$ with leaves on the same three distinct components of
$\ti L$ does not change the surgery equivalence class of $L$.
First of all we can consider the case where $G$ is an inverse of $G'$ in the
sense of \cite[Theorem 3.2]{GGP}. In this case a surgery on $G$ and $G'$ does
not change $L$ at all. For any other $G$ we can assume that there is a homotopy
in $N$ from $G$ to an inverse of $G'$  which is stationary on the leaves of 
$G$.
Such a homotopy is a sequence of isotopies in $N-L$ together with (i) crossings
of an edge of $G$ and a component of $L$ and (ii) crossings of an edge of $G$
with a leaf of $G$. It suffices to show that these two types of crossings do 
not change the surgery equivalence class of the $G$ surgery on $L$.

For (i) the effect of this crossing on surgery of $L$ is pictured in Figure
\ref{fig.delcross}. 

\begin{figure}[ht!]
\psdraw{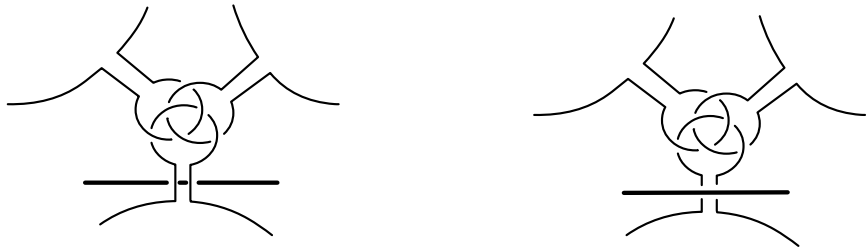}{4.5in} 
\caption{A link component crossing an edge of a Y--graph}\lbl{fig.delcross}
\end{figure}

The surgery equivalence is given by Lemma \ref{lem.conc} and
the double crossing change in Figure \ref{fig.surger1}. 

\begin{figure}[hb!]
\psdraw{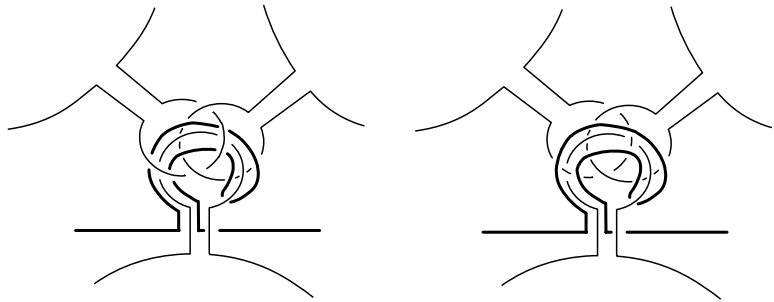}{4in} 
\caption{A double crossing change which implements the crossing of Figure
\ref{fig.delcross}}\lbl{fig.surger1}
\end{figure}

For (ii) we invoke the
following:
\begin{lemma}
Suppose $G$ is the union of two \ylink s $G_1 ,G_2$  in the complement of a 
link
$L$, whose leaves are meridians of $L$, and $G'$ is obtained from $G$ by a
single crossing change of a leaf of $G_1$ with a leaf of $G_2$. Then the link
produced by surgery on $L$ using $G$ is surgery equivalent to the link produced
by surgery on $G'$.
\end{lemma}
\begin{proof} By \cite[Theorem 2.3]{GGP} and $Y_4$ moves of \cite{GGP}, surgery
on $G'$ is the same as surgery on $G$ together with surgery on a clover of
degree $2$ with the shape of $\ppsdraw{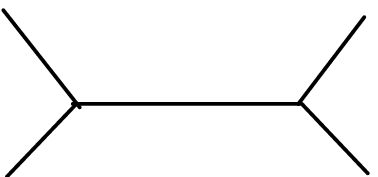}{.4in}$.
Thus we need to see that surgery on such clovers does not change the surgery
equivalence class of $L$. The effect of surgery on such a clover is shown in
Figure
\ref{fig.hlink}. 

\begin{figure}[htpb!]
\psdraw{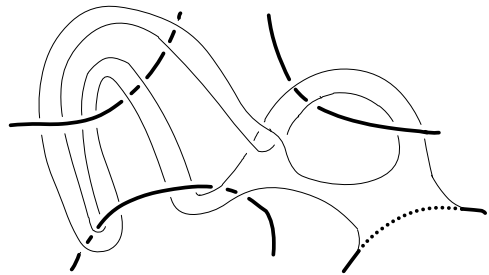}{2.5in} 
\caption{The effect of surgery on a clover of degree $2$}\lbl{fig.hlink}
\end{figure}

A double crossing change which will undo this surgery is illustrated in Figure
\ref{fig.surger2}.\end{proof}

\begin{figure}[ht!]
\psdraw{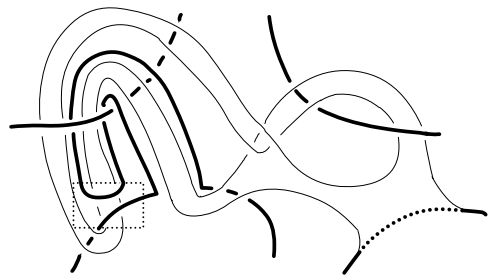}{2.5in} 
\caption{The doublecrossing change inside the box will undo the
surgery.}\lbl{fig.surger2}
\end{figure}

This completes the proof of Theorem \ref{th.surg}.
\end{proof}

We also prove that:

\begin{theorem}
\lbl{thm.conc}
Concordant links are surgery equivalent.
\end{theorem}

\begin{proof}
Concordance, just as in the classical case, is generated by the following 
{\em ribbon move} $L\to L'$. Given a $\pi$--AS link $L\sub N$, consider
also a finite number of disks $\{ D_i\}$ in $N$, disjoint from each other 
and $L$. For each $D_i$ choose a band $B_i$ connecting $D_i$ to a component,
which we denote $L_i$, of $L$. The band cannot intersect $L$ or any $\bd D_j$
except at its ends. Then $L'$ is defined to be the band-sum of $L$ with 
$\bigcup\bd D_i$.

Now choose some place where a band $B_i$ penetrates a disk $D_j$. Choose
a path $\g$ from $L_j$ to nearby the penetration so that the closed path 
consisting of $\g$ followed by the path from  $D_j$ along $B_j$ and back 
along $L_j$ to the starting point of $\g$ is null-homotopic.
Now thicken $\g$ to a band (or finger) and apply Lemma \ref{lem.conc}
as follows:

\psdraw{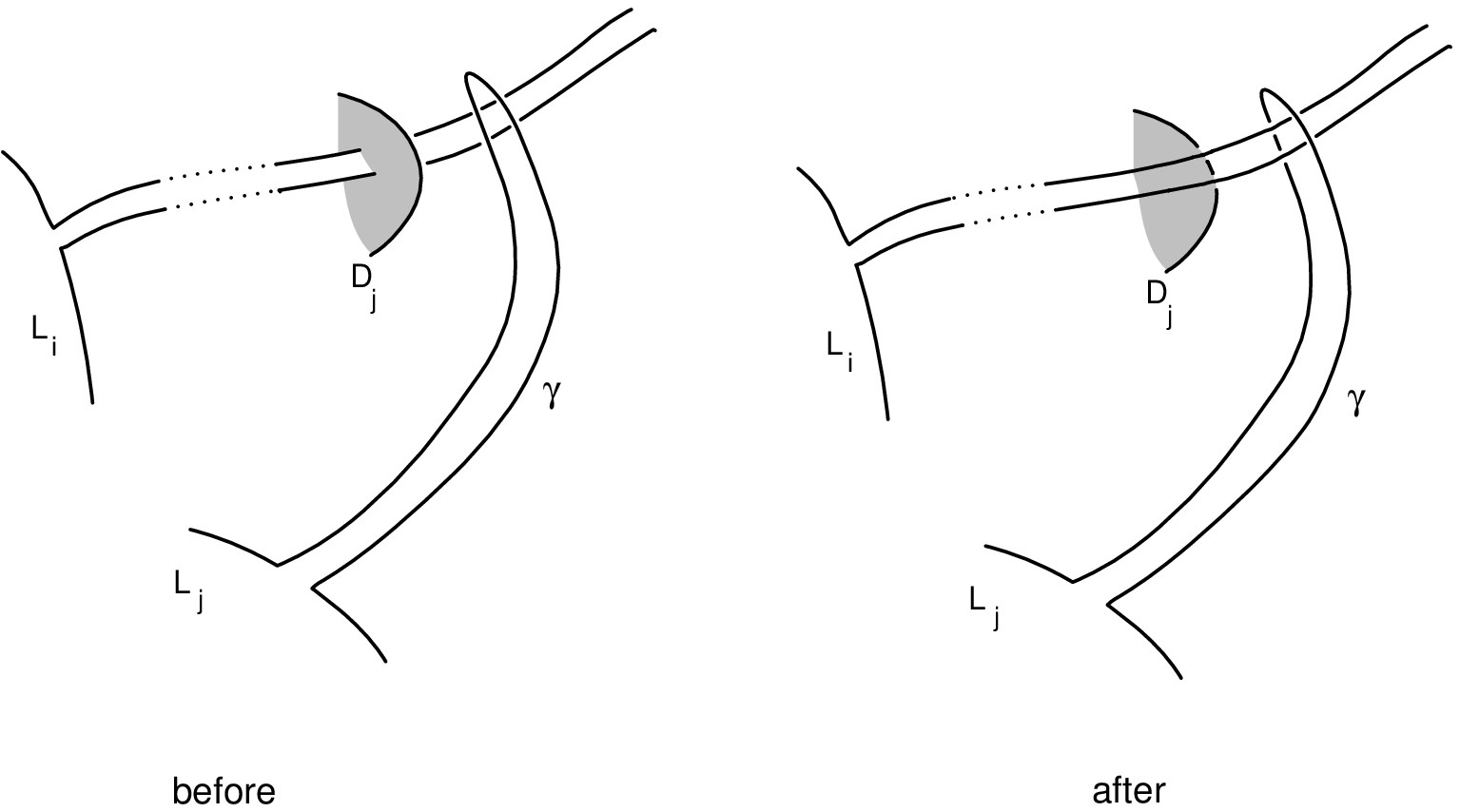}{4in}

\noi This removes the penetration. Eventually we can remove all the 
penetrations and the resulting link will be isotopic to $L$.
\end{proof}

\subsection{An alternative study of $\Gas {}$.}
\lbl{sub.nequiv}

In this section we mention, in brief, an alternative study of $\GY {}$
using our results on surgery equivalence and the group $\A (\pi )$ from
Section \ref{sec.gradedq}. For $N=S^3$ this coincides with the
approach to \fti s introduced by Ohtsuki \cite{Oh}, and studied in \cite{GO}.

Recall the map $\Las(N)\to\ti\K(N)$ defined by doing surgery on a unit-framed 
$\pi$--AS link in $N$. If $L$ and $L'$ are surgery equivalent $n$ component 
$\pi$--AS links in $N$, then $[N,L]=[N,L'] \bmod \F_{n+1}(N)$.
Thus $\G(N)$ is a quotient of the free abelian group on the set of
surgery equivalence classes of $\pi$--AS links. 
 For each
generator $(G,\a)$ of $\A(\pi )$ we can
define a $\pi$--AS link as follows. First construct a link in a $3$--ball $B\sub
N$ associated to $G$, using the construction of Ohtsuki \cite{Oh}, by banding
together copies of the Borromean rings, one Borromean rings for each vertex and
one band for each edge. Now for each edge $e$ of $G$ pull the band corresponding
to $e$ around a loop representing $\a(e)$, using the orientation of $e$ to
direct it. See Figure \ref{fig.edge} below. 
\begin{figure}[ht!]
\psdraw{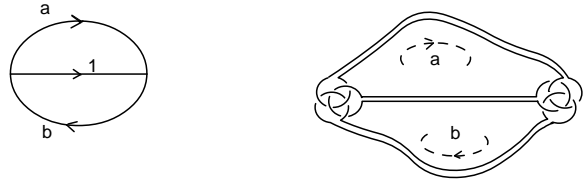}{3in}
\caption{An edge and vertex oriented diagram and the associated $\pi$--AS
link}\lbl{fig.edge}
\end{figure}
The surgery equivalence class of this link $L(G)$ is well-defined since we only
have to worry about band-crossings which are covered by Lemma \ref{lem.conc}.
Note that if $G$ has degree $2n$, then $L(G)$ has $3n$ components. 
Using a local {\em $3$--band relation} of \cite[Lemma 4.1]{Oh} for an 
arbitrary manifold $N$ and the above discussion, we obtain an onto
map from the abelian group generated by pairs $(G,\a)$ of $\text{degree}(G)=
2n$ to $\Gas {3n}\oZ$. The work of \cite{GO}, formulated for arbitrary
manifolds rather than $S^3$, implies that the $\AS$ and $\IHX$ relations
of Figure \ref{relations} are satisfied, thus obtaining an onto map
$\A_{2n}(\pi)\twoheadrightarrow \Gas {3n}\oZ$ for every integer $n$.
It is an easy exercise in Kirby calculus to show that the above maps
fit in a commutative diagram
$$
\begin{diagram}
\node{\A_{2n}(\pi
)}\arrow{e,A}\arrow{s,r,<>}{\iso}\node{\G_{3n}(N)\oZ}\arrow{s,r,<>}{\iso}\\
\node{\A'_{2n}(\pi )}\arrow{e,A}\node{\G^Y_{2n}(N)\oZ}
\end{diagram}
$$

{\bf Acknowledgement}\qua The  authors were partially supported by NSF grants
        DMS-98-00703 and DMS-99-71802 respectively, and by an Israel-US
BSF grant.


\begin{thebibliography}

\let\oldbibitem\bibitem
\def\bibitem#1]#2{\oldbibitem{#2}}


\bibitem[FR]{FR} {\bf R Fenn}, {\bf C Rourke}, 
        {\em On Kirby's calculus of links}, 
        Topology, {18} (1979) 1--15

\bibitem[FQ]{FQ} {\bf M Freedman}, {\bf F Quinn}, 
        {\it Topology of 4--manifolds}, 
        Princeton University Press, Princeton, NJ (1990)
        
\bibitem[GGP]{GGP} {\bf S Garoufalidis},  {\bf M Goussarov}, {\bf M Polyak},
        {\em Calculus of \clover s and finite type invariants of
        3--manifolds}, 
        Geometry and Topology, {5} (2001) 75--108

\bibitem[GL]{GL} {\bf S Garoufalidis}, {\bf J Levine}, 
        {\em Finite-type invariants of 3--manifolds II},  
        Math. Annalen {306} (1996), 691--718
        
\bibitem[GO]{GO} {\bf S Garoufalidis}, {\bf T Ohtsuki},
        {\em On finite type 3--manifold invariants III: manifold weight 
        systems}, 
        Topology, {37} (1998) 227--244
        
\bibitem[Gu]{Gu} {\bf M Goussarov},
        {\em Finite type invariants and $n$--equivalence of
        3--manifolds}, 
        CR Acad. Sci. Paris Ser. I Math. {329} (1999) 517--522

\bibitem[H]{H} {\bf K Habiro}, 
        {\em Claspers and finite type invariants of links},
        Geometry and Topology, {4} (2000) 1--83

\bibitem[KT]{KT} {\bf R Kirby}, {\bf L\,R Taylor},
        {\em A survey of 4--manifolds though the eyes of surgery},
        from: ``Surveys on surgery theory'', 
        vol. 1, ed. S Cappell, A Ranicki, J.                             
        Rosenberg, Annals of Math. Studies 145, Princeton U. Press (2000)

\bibitem[Le]{Le} {\bf J Levine},
        {\em Surgery equivalence of links}, 
        Topology, {26} (1987) 45--61

\bibitem[LMO]{LMO}  {\bf T\,T\,Q  Le}, {\bf J Murakami}, {\bf T Ohtsuki},
        {\em A universal quantum invariant of 3--manifolds},
        Topology, {37} (1998) 539--574 

\bibitem[Ma]{Ma} {\bf S\,V Matveev},
        {\em Generalized surgery of three-dimensional manifolds and
        representations of homology spheres},
        Math. Notices Acad. Sci. USSR,
        {42:2} (1987) 651--656

\bibitem[MN]{MN} {\bf H Murakami}, {\bf Y Nakanishi},
        {\em On a certain move generating link homology},
        Math. Annalen, {284} (1989) 75--89

\bibitem[Oh]{Oh} {\bf T Ohtsuki},
        {\em Finite type invariants of integral homology 3--spheres}, 
        J. Knot Theory and its  Ramifications, {5} (1996) 101--115 

\bibitem[R]{R} {\bf A Ranicki},
        {\em The algebraic theory of surgery I. Foundations}, 
        Proc. London Math. Soc. {40} (1980) 87--192 
        
\bibitem[Wa]{Wa} {\bf C\,T\,C Wall},
        {\it Surgery on manifolds}, 
        Academic Press, London (1970)

        
\end{thebibliography}
\end{document}
\endinput